\documentclass[oneside]{amsart}
\usepackage[colorlinks]{hyperref} 
\usepackage[backend=biber, style=alphabetic, sorting=nyt, maxbibnames=99, maxalphanames=4, minalphanames=3, doi=false, isbn=false, url=false, eprint=false]{biblatex} 
\usepackage{xurl}
\hypersetup{breaklinks=true}

\DeclareFieldFormat*{title}{\mkbibemph{#1}}
\DeclareFieldFormat*{journaltitle}{#1}
\DeclareFieldFormat*{booktitle}{#1}

\renewbibmacro{in:}{%
  \ifentrytype{article}
    {}
    {\printtext{\bibstring{in}\intitlepunct}}}

\DeclareFieldFormat[article,periodical]{volume}{\mkbibbold{#1}}
\DeclareFieldFormat[article,periodical]{number}{\bibstring{number}~#1}
\renewbibmacro*{journal+issuetitle}{%
  \usebibmacro{journal}%
  \setunit*{\addspace}%
  \iffieldundef{series}
    {}
    {\newunit
     \printfield{series}%
     \setunit{\addspace}}%
  \printfield{volume}%
  \setunit{\addspace}%
  \usebibmacro{issue+date}%
  \setunit{\addcolon\space}%
  \usebibmacro{issue}%
  \setunit{\addcomma\space}%
  \printfield{number}%
  \setunit{\addcomma\space}%
  \printfield{eid}%
  \newunit}

\DeclareFieldFormat[article,periodical]{pages}{#1}

\DeclareFieldFormat{mrnumber}{%
  MR\addcolon\space
  \ifhyperref
    {\href{http://www.ams.org/mathscinet-getitem?mr=#1}{\nolinkurl{#1}}}
    {\nolinkurl{#1}}}

\renewbibmacro*{doi+eprint+url}{%
  \iftoggle{bbx:doi}
    {\printfield{doi}}
    {}%
  \newunit\newblock
  \printfield{mrnumber}%
  \newunit\newblock
  \iftoggle{bbx:eprint}
    {\usebibmacro{eprint}}
    {}%
  \newunit\newblock
  \iftoggle{bbx:url}
    {\usebibmacro{url+urldate}}
    {}}

\DeclareFieldFormat{postnote}{\mknormrange{#1}}
\DeclareFieldFormat{volcitepages}{\mknormrange{#1}}
\DeclareFieldFormat{multipostnote}{\mknormrange{#1}}

\usepackage{booktabs} 

\let\originallhook=\lhook

\usepackage{amsthm,amssymb,amsmath,mathrsfs,stmaryrd}  
\usepackage[LGR,T1]{fontenc}  
\usepackage{lmodern}
\usepackage{xcolor}
\usepackage[all,matrix,color,arrow,pdf]{xy}  
\usepackage{graphicx}
\usepackage{mathtools}
\usepackage{geometry}
\usepackage{tikz}  
\usetikzlibrary{calc, intersections}
\usepackage{enumerate}
\usepackage{pdfpages}
\usepackage[utf8]{inputenc}
\usepackage{longtable}
\usepackage{setspace}
\usepackage{hyperref}
\usepackage{eucal} 
\usepackage{tikz-cd}
\usepackage{slashed}
\usepackage{relsize}
\usepackage[bbgreekl]{mathbbol}
\usepackage{amsfonts}
\hypersetup{linkcolor=[RGB]{144,0,32}, citecolor=[RGB]{33,66,171}}

\DeclareSymbolFontAlphabet{\mathbb}{AMSb}
\DeclareSymbolFontAlphabet{\mathbbl}{bbold}
\newcommand{\prism}{\mathbbl{\Delta}}

\let\lhook=\originallhook

\newcommand\sbullet[1][.5]{\mathbin{\vcenter{\hbox{\scalebox{#1}{$\bullet$}}}}} 

\DeclareMathAlphabet{\mathpzc}{OT1}{pzc}{m}{it}

\def\ker{\operatorname{ker}}

\def\Mor{\operatorname{Mor}}
\def\Map{\operatorname{Map}}
\def\map{\underline{\operatorname{Map}}}

\def\Sp{\operatorname{Sp}}
\def\Spec{\operatorname{Spec}}

\def\Spf{\operatorname{Spf}}

\def\tr{\operatorname{tr}}

\def\K{\mathcal{K}}

\def\Dh{\widehat{\mathcal{D}}}

\def\DFh{\widehat{\mathcal{DF}}}

\def\op{\mathrm{op}}

\def\N{\operatorname{N}}

\def\colim{\operatorname{colim}}

\def\Ab{\text{Ab}}
\def\Ring1{\text{Ring}1}
\def\CRing1{\text{CRing}1}
\def\Mod{\text{Mod}}

\def\CAlg{\text{CAlg}}

\def\an{\text{an}}

\def\bbE{\mathbb{E}}
\def\bbF{\mathbb{F}}
\def\bbG{\mathbb{G}}

\def\bbZ{\mathbb{Z}}

\def\bbS{\mathbb{S}}

\def\+1{\xrightarrow{+1}}

\def\dt{\otimes^{\mathrm{L}}}

\def\RG{\text{R}\Gamma}
\def\Rlim{\text{R}\lim}

\def\lsh{_{!}}

\def\Sym{\operatorname{Sym}}

\def\aff{\mathrm{aff}}
\def\zar{\mathrm{Zar}}

\def\et{\mathrm{\'et}}
\def\proet{\mathrm{pro\'et}}

\def\crys{\text{crys}}
\def\dR{\text{dR}}

\def\gr{\operatorname{gr}}

\def\isomto{\overset{\sim}{\to}}

\def\Sp{\text{Sp}}

\def\N{\text{N}}
\def\Fun{\text{Fun}}
\def\fib{\text{fib}}
\def\cof{\text{cof}}

\def\Pro{\text{Pro}}

\def\Cpl{\operatorname{Cpl}}

\def\Shv{\operatorname{Shv}}

\def\cn{\mathrm{cn}}

\def\Poly{\mathrm{CAlg}^{\text{poly}}}

\def\poly{\mathrm{poly}}

\def\THH{\mathrm{THH}}
\def\TC{\mathrm{TC}}

\def\TR{\mathrm{TR}}
\def\tr{\mathrm{tr}}

\def\TP{\mathrm{TP}}

\def\K{\mathrm{K}}

\def\Fil{\mathrm{Fil}}
\def\fil{\mathrm{fil}}
\def\disc{\mathrm{disc}}

\def\Nyg{\mathcal{N}}
\def\QSyn{\mathrm{QSyn}}

\def\dR{\mathrm{dR}}
\def\dRh{\widehat{\mathrm{dR}}}
\def\conj{\mathrm{conj}}

\theoremstyle{definition}
\newtheorem{example}{Example}[section]  
\newtheorem{definition}[example]{Definition}

\newtheorem{proposition}[example]{Proposition}
\newtheorem{lemma}[example]{Lemma}

\newtheorem{theorem}[example]{Theorem}
\newtheorem{corollary}[example]{Corollary}

\newtheorem{remarkn}[example]{Remark} 

\newtheorem*{ack}{Acknowledgements}
\newtheorem*{conv}{Conventions}

\theoremstyle{remark}
\newtheorem*{remark}{Remark}

\begin{document}
\title{Topological cyclic homology of Cartier smooth rings}
\author{Hyungseop Kim}
\address{Department of Mathematics, University of Toronto, Bahen Centre, 40 St. George St., Toronto, ON, Canada, M5S 2E4}
\email{khsato@math.utoronto.ca}

\begin{abstract}
We study algebraic K-theory, syntomic cohomology, and prismatic cohomology of Cartier smooth rings. As an application, we provide an alternative proof of Kelly-Morrow's generalization of the Geisser-Levine theorem computing $p$-adic algebraic K-theory of Cartier smooth local rings; our approach relies on the description of topological cyclic homology through the motivic filtration. 
\end{abstract}

\maketitle

\tableofcontents
 
\section{Introduction}
In this article, we study algebraic K-theory, syntomic cohomology, and prismatic cohomology complexes of Cartier smooth rings. In particular, we explain an alternative approach to Kelly-Morrow's generalization of the Geisser-Levine theorem to the Cartier smooth case through the motivic filtration on topological cyclic homology. Let us first recall what the theorem states. Fix a prime $p$. Thanks to Quillen's classical computation \cite{quil} of algebraic K-theory of finite fields, we have the following concrete description of algebraic K-theory groups in all degrees:
\begin{equation*}
\K(\mathbb{F}_{p}) : \begin{cases} 
                   \K_{<0}(\mathbb{F}_{p}) = 0,\\
                   \K_{0}(\mathbb{F}_{p}) = \mathbb{Z},\\
                   \K_{2i}(\mathbb{F}_{p}) = 0,\\
                   \K_{2i-1}(\mathbb{F}_{p}) = \mathbb{Z}/(p^{i}-1).
                    \end{cases}
\end{equation*}
Thus, we have vanishing of positive even degrees and $p$-torson freeness in all degrees. In particular, $\K(\mathbb{F}_{p},\mathbb{Z}_{p}) := \K(\mathbb{F}_{p})^{\wedge_{p}} \simeq H\mathbb{Z}_{p}$ as an $\bbE_{\infty}$-ring; this can be seen as a baby example of chromatic redshift for algebraic K-theory. Note that while the input $\bbF_{p}$ is in characteristic $p$, after $p$-completion the output algebraic K-theory looks like a lifting of a characteristic $p$ object to the mixed characteristic $(0,p)$. One might hope that, although computing (integral) algebraic K-groups would generally be difficult, describing $p$-adic algebraic K-theory of characteristic $p$ rings akin to the consequence $\K(\bbF_{p})^{\wedge_{p}}\simeq H\bbZ_{p}$ of Quillen's theorem might be more approachable. In fact, Geisser and Levine showed the following generalization of the aforementioned $p$-adic consequence of Quillen's computation in the smooth case. 

\begin{theorem}[\cite{gl}]\label{thm:glintro}
Let $R$ be a local smooth algebra over a perfect field $k$ of characteristic $p$. Then, for each $i\geq 0$ and $r\geq 1$, the followings hold:\\
(1) $\K_{i}(R)$ is $p$-torsion free. \\
(2) $\K_{i}(R)/p^{r}\simeq W_{r}\Omega^{i}_{R,\log}$. 
\end{theorem}

Here, $W_{r}\Omega^{i}_{R,\log}$ stands for the logarithmic part of the $i$-th de Rham-Witt group of length $r$ \cite{ill}. Under the $p$-torsion freeness condition of (1), the condition (2) amounts to the equivalence $\K_{i}(R,\bbZ_{p})\simeq W\Omega^{i}_{R,\log}$, where $\K_{i}(R,\bbZ_{p}) = \pi_{i}(\K(R)^{\wedge_{p}})$ as usual. As the de Rham-Witt complex $W\Omega^{\ast}_{R}$ provides a lifting of the de Rham complex $\Omega^{\ast}_{R}$ of $R$ to mixed characteristic $(0,p)$, one observes that the pattern observed in the case of $R=\bbF_{p}=k$ remains valid. 

\begin{example}\label{ex:introFp}
Let $R = \bbF_{p}$. In this case, the de Rham-Witt complex is given by $W\Omega^{\ast}_{\bbF_{p}}\simeq \widehat{\Omega}^{\ast}_{\mathbb{Z}_{p}}$, the $p$-completed de Rham complex of $W(R) = \bbZ_{p}$. Moreover, since $L_{\bbZ_{p}/\bbZ}^{\wedge_{p}}\simeq 0$, we know $W\Omega^{\ast}_{\bbF_{p}}\simeq \widehat{\Omega}^{\ast}_{\mathbb{Z}_{p}}\simeq \bbZ_{p}$. Note that the Frobenius map on $W\Omega^{\ast}_{\mathbb{F}_{p}}$ is $F=\varphi = id_{\bbZ_{p}}$, and hence we have 
\begin{equation*}
W\Omega^{i}_{\bbF_{p},\log} \simeq \begin{cases} 
           \bbZ_{p}, & i=0,~\text{and} \\
           0, & i>0.
           \end{cases}
\end{equation*}
Through the Geisser-Levine theorem, this recovers $\K(\bbF_{p},\bbZ_{p})\simeq H\bbZ_{p}$. Cf. Lemma \ref{lem:perfdW} for relevant details and generalizations. 
\end{example} 

\begin{remarkn}
Note that in the Geisser-Levine description of algebraic K-groups, both sides of the equivalence commute with filtered colimits of rings. Thus, Theorem \ref{thm:glintro} remains valid for local ind-smooth $k$-algebras. In particular, Theorem \ref{thm:glintro} applies to regular local rings in characteristic $p$, since these rings are ind-smooth over $\bbF_{p}$ by N\'eron-Popescu desingularization. 
\end{remarkn}

Another class of rings which contain characteristic $p$ finite fields is that of perfect rings. As arguments of Example \ref{ex:introFp} look compatible with perfect rings (more or less trivially, up to verifying that the de Rham-Witt complex should be the ring of $p$-typical Witt vectors), or more interestingly with smooth algebras over perfect rings (again once de Rham-Witt complex of such rings are properly defined), it is tempting to speculate that Geisser-Levine description should remain true for smooth algebras over perfect rings. In \cite{km}, Kelly and Morrow proved that that is indeed the case. In fact, they introduced the notion of Cartier smoothness for characteristic $p$ rings, and proved that Geisser-Levine description holds for local Cartier smooth rings. 

\begin{theorem}[\cite{km}]\label{thm:kmintro}
Let $S$ be a local Cartier smooth algebra over $\bbF_{p}$. Then, for each $i\geq 0$ and $r\geq 1$, the followings hold:\\
(1) $\K_{i}(S)$ is $p$-torsion free. \\
(2) $\K_{i}(S)/p^{r}\simeq W_{r}\Omega^{i}_{S,\log}$. 
\end{theorem}

Contemporaneous with \cite{km}, Cartier smoothness was also studied in \cite{vorst}. Archetypal examples of Cartier smooth rings include smooth algebras over perfect rings and valuation rings containing $\bbF_{p}$. We will review the notion of Cartier smoothness and some of its characterizations in sections \ref{sec:cartsm} and \ref{sec:cartsmderW}. \\
\indent In \cite{km}, Theorem \ref{thm:kmintro} was proved by reduction to Theorem \ref{thm:glintro} of Geisser-Levine through trace methods. More precisely, they first handle homotopy groups $\pi_{\ast}(\TC(S)/p^{r})$ via computations of topological restriction theory groups $\pi_{\ast}(\mathrm{TR}(S)/p^{r})$ and Frobenius on them, using the classical description of $\TC/p^{r}$ as a Frobenius fixed point of $\TR/p^{r}$. Then, they use CMM rigidity \cite{cmm} in trace methods to ensure that the difference between algebraic K-theory groups and logarithmic de Rham-Witt groups for local Cartier smooth rings, where the latter is related to $\TC/p^{r}$ by the previous computation, is exactly the same as the difference between them in the case of ind-smooth algebras, which is nothing but zero due to the Geisser-Levine theorem.  \\
\indent Our approach to Theorem \ref{thm:kmintro} follows the same strategy of reducing the problem to the case of ind-smooth algebras through trace methods. However, our way of describing $\TC$ and treating the entire situation $p$-completely requires us to employ different techniques and offers us additional useful information complementary to \cite{km}. First, we approach topological cyclic homology of Cartier smooth rings $\TC(S,\bbZ_{p})$ through its motivic filtration \cite{bms2}\footnote{We note that it was already indicated in \cite{km} that it would be possible to use the motivic filtration to reformulate some of their computations.}. More precisely, we compute syntomic cohomology complexes $\bbZ_{p}(n)(S)$ of Cartier smooth rings $S$, which are associated graded pieces of the motivic filtration on $\TC(S,\bbZ_{p})$ and hence act as $p$-adic (pro-)\'etale motivic cohomology complexes of $\Spec S$. It turns out syntomic cohomology complexes are precisely given by logarithmic de Rham-Witt sheaves up to shift, and this gives a description of $\TC(-,\bbZ_{p})$ in terms of a pro-\'etale sheaf on $\Spec S$:

\begin{proposition}[Corollary \ref{cor:bms2log2} and Corollary \ref{cor:cartsmTCpadic}] \label{prop:introsyntomic}
Let $X$ be a Cartier smooth $\bbF_{p}$-scheme. Then, $\bbZ_{p}(n)_{X_{\proet}}\simeq W\Omega^{n}_{X,\log}[-n]$ and $\pi_{n}^{\proet}\TC(-,\bbZ_{p})\simeq W\Omega^{n}_{X,\log}$ for all $n\geq0$ as pro-\'etale sheaves on $X$. 
\end{proposition} 

We note that the computation of the (mod-$p^{r}$) syntomic cohomology complexes for Cartier smooth rings has already appeared in \cite[Prop. 5.1 (ii)]{lm}. This can be viewed as a special case of \cite[Th. 1.8]{bm} for the class of F-smooth rings, which provides a mixed characteristic generalization of Cartier smoothness; cf. Remark \ref{rem:F-smooth} below. \\
\indent Next, we use pro-\'etale descent spectral sequence to obtain natural exact sequences relating homotopy groups of $p$-adic $\TC$ and syntomic cohomology complexes, cf. Proposition \ref{prop:proetpostnikov}. Finally, we use CMM rigidity as well as AMMN rigidity property \cite{ammn} for syntomic cohomology complexes to compute $\pi_{i}\K_{\geq 0}(S,\bbZ_{p})$ in Theorem \ref{th:cartsmglpadic}, which implies Theorem \ref{thm:kmintro}:

\begin{theorem}[Theorem \ref{th:cartsmglpadic}]
Let $S$ be a Cartier smooth local ring. Then, there is a natural isomorphism $\K_{i}(S,\bbZ_{p})\simeq W\Omega^{i}_{S,\log}$ induced by the cyclotomic trace map for all $i\geq 0$. 
\end{theorem}

The key point of the computation of the syntomic cohomology complexes of Cartier smooth rings given here is that the prismatic cohomology complexes $\prism_{S}$ of those rings are precisely given by their de Rham-Witt complexes together with the data of Nygaard filtrations: 

\begin{proposition}[Proposition \ref{prop:cartsmLdW} and Corollary \ref{cor:cartsmprism}] \label{prop:introprism}
Let $X$ be a Cartier smooth $\bbF_{p}$-scheme. Then, there are natural equivalences $\Fil^{\sbullet}_{\N}\RG_{\prism}(X)\simeq \RG(X,\Nyg^{\geq\sbullet}\mathcal{W}\Omega_{X})\simeq \RG(X,\Nyg^{\geq\sbullet}W\Omega_{X})$ in $\DFh(\bbZ_{p})$. 
\end{proposition} 

Here, $\mathcal{W}\Omega_{X}$ stands for the saturated de Rham-Witt complex of $X$ as studied in \cite{blm}. Thus, we in particular know that for Cartier smooth rings, saturated de Rham-Witt complex $\mathcal{W}\Omega_{S}$ agrees with the underived de Rham-Witt complex $W\Omega_{S}$ as well as with the prismatic cohomology complex $\prism_{S}$. We will remark on this fact and relevant consequences involving $p$-torsion free lifts of Cartier smooth rings (which are instances of mixed characteristic F-smooth rings) in section \ref{sec:cartsmprismetc}. We also note that for F-smooth rings over perfectoid rings, \cite[Th. C]{bouis} provides computations concerning their prismatic cohomology complexes. In this article, we provide a unified approach towards Proposition \ref{prop:introsyntomic} and Proposition \ref{prop:introprism} involving the saturated de Rham-Witt complexes of \cite{blm} for Cartier smooth rings. Based on this, we present a more structured approach to reproving the Geisser-Levine theorem for Cartier smooth local rings \cite{km} using the Nikolaus-Scholze description \cite{ns} of $\TC$ and its motivic filtration \cite{bms2}.

\begin{conv}
We fix a prime number $p$. For each commutative ring $R$, we denote $\mathcal{D}(R)\simeq \Mod_{R}$ for its derived $\infty$-category. When $R$ is $p$-complete, we denote $\widehat{\mathcal{D}}(R)\simeq\Mod_{R}^{\Cpl(p)}$ for its $p$-complete derived $\infty$-category, and similarly $\widehat{\mathcal{DF}}(R)\simeq\widehat{\mathcal{D}}(R)^{\fil}$ for the filtered $p$-complete derived $\infty$-category of $R$. 
\end{conv}

\begin{ack}
We would like to thank Matthew Morrow for helpful communications. This article is based on a part of the author's doctoral thesis \cite{kimthesis} with minor expositional changes, and the author would like to thank Kay R\"ulling for helpful comments as the External Appraiser and Michael Groechenig for support. The author was supported by the CMK foundation. 
\end{ack}

\section{Cartier smoothness}\label{sec:cartsm}

In this section, we briefly recall the notion of Cartier smoothness for characteristic $p$ rings as introduced in \cite{km} (and also studied in \cite{vorst}). Recall that for a commutative ring $S$ of characteristic $p$, the inverse Cartier map $C^{-1}:\Omega^{\sbullet}_{S^{(1)}/\bbF_{p}}\to H^{\sbullet}(\Omega^{\ast}_{S/\bbF_{p}})$ is a $\bbZ_{\geq0}$-graded $S^{(1)}$-algebra map determined by the requirements $C^{-1}(f) = f^{p}$ and $C^{-1}(dg) = [g^{p-1}dg]\in H^{1}(\Omega^{\ast}_{S/\bbF_{p}})$ for $f,g\in S$. 

\begin{definition}(\cite{km})
A commutative $\bbF_{p}$-algebra $S$ is called \emph{Cartier smooth} if its cotangent complex $L_{S/\bbF_{p}} = L\Omega^{1}_{S/\bbF_{p}}$ is a flat (ordinary) $S$-module and the inverse Cartier map $C^{-1}:\Omega^{\sbullet}_{S^{(1)}/\bbF_{p}}\to H^{\sbullet}(\Omega^{\ast}_{S/\bbF_{p}})$ is an isomorphism of $\bbZ_{\geq0}$-graded rings. An $\bbF_{p}$-scheme $X$ is \emph{Cartier smooth} if for every affine open subset of $X$, its ring of functions is Cartier smooth. 
\end{definition}

In particular, the cotangent complex $L_{S/\bbF_{p}}$ of Cartier smooth $S$ has Tor-amplitude in cohomological degrees $[-1,0]$, i.e., $S$ is quasisyntomic. 

\begin{example}
The following class of rings are Cartier smooth: \\
(1) (ind-)smooth algebras over a perfect field $k$ of characteristic $p$. \\
(2) perfect rings. \\
(3) (ind-)smooth algebras over a perfect ring. \\
(4) filtered colimits and localizations of Cartier smooth rings. \\
(5) (N\'eron-Popescu) regular Noetherian rings over $\bbF_{p}$. \\
(6) (Gabber \cite[App. A]{vorst}) valuation rings over $\bbF_{p}$. 
\end{example}

\begin{remarkn}
(1) A Noetherian $\bbF_{p}$-algebra is Cartier smooth if and only if it is regular, cf. \cite[Th. 4.15]{bm} (in the Cartier smooth case, the main technical ingredient is in \cite[Prop. 4.19]{bm}). Thus, Cartier smoothness can be regarded as a non-Noetherian generalization of regularity in characteristic $p$.\\
(2) Despite (1), note that Cartier smooth rings in general have non-vanishing negative K-groups, e.g., perfect rings might have nonzero negative K-groups. In the setting of Geisser-Levine Theorem \ref{thm:glintro}, ind-smooth $k$-algebras of course have vanishing negative K-groups, so describing the connective part was already enough. Theorem \ref{thm:kmintro} gives a description of the connective $p$-adic algebraic K-theory for Cartier smooth rings, but says nothing about the negative K-group. \\
(3) Both \cite{km} and our approach to Theorem \ref{thm:kmintro} relies on the trace method in the $p$-adic setting, which by its very nature can only control the connective part of algebraic K-theory. This might be the main reason for the limitation described in (2).  
\end{remarkn}

\begin{remarkn} \label{rem:F-smooth}
For general ($p$-adic) situation, \cite{bm} introduced the notion of F-smoothness for any $p$-quasisyntomic rings in order to capture an absolute version of $p$-complete smoothness\footnote{Roughly speaking, the condition of F-smoothness is designed to enforce the Segal conjecture true under the presense of finiteness of conjugate filtration.}; over fixed perfectoid base rings, this was also studied in \cite{bouis}. This notion precisely generalizes the notion of Cartier smoothness in characteristic $p$; a quasisyntomic $\bbF_{p}$-algebra is Cartier smooth if and only if it is F-smooth \cite[Prop. 4.14]{bm}. Moreover, in \emph{loc. cit.} they computed (mod-$p^{r}$) syntomic complexes of $p$-torsion free F-smooth rings in terms of truncated nearby cycles (which, as usual, is the \'etale Tate twist on the locus where $p$ acts invertibly) modified by the image of the symbol map (which is the only visible part on the vanishing locus of $p$, exactly as in Corollary \ref{cor:bms2log2}) at the top cohomological degree \cite[Th. 1.8]{bm}; this is a $p$-adic generalization of our computation of syntomic complexes in Corollary \ref{cor:bms2log2} for the characteristic $p$ case, albeit exact methods of computations are different. 
\end{remarkn}

\section{De Rham-Witt complex of Cartier smooth rings}\label{sec:cartsmderW}

In this section, we recall the notion of Cartier smoothness and study de Rham-Witt complex of Cartier smooth rings. The slogan is that every de Rham-type invariant is non-derived, i.e., derived objects agree with their underived analogues, for Cartier smooth rings. De Rham-Witt complexes \cite{ill} of Deligne-Illusie are objects $W\Omega^{\ast}_{X}$ naturally constructed for each $\bbF_{p}$-scheme $X$ which, in the smooth case, provide an explicit complex of sheaves whose (hyper-)cohomology computes crystalline cohomology of Berthelot-Grothendieck for each smooth $\bbF_{p}$-scheme. Similarly, one can consider the derived de Rham-Witt complex $LW\Omega_{X}$ which computes derived crystalline cohomology of a (derived) scheme $X$ over $\bbF_{p}$. In order to remedy some computational and conceptual difficulties involved in \cite{ill}, an alternative construction of de Rham-Witt complexes was proposed by \cite{blm}; for each $\bbF_{p}$-scheme $X$, they attach the saturated de Rham-Witt complex $\mathcal{W}\Omega^{\ast}_{X}$ which agrees with Illusie's $W\Omega^{\ast}_{X}$ whenever $X$ is smooth. By construction saturated de Rham-Witt complexes have arguably simpler universal property and enable one to detour technical complications from the use of pro-complexes (which \cite{ill} cannot avoid) at the cost of being difficult to carry out explicit computations outside of the smooth case and to describe (co)limits in general due to saturation. Nevertheless, we will show that for Cartier smooth rings (and hence for Cartier smooth schemes), all three de Rham-Witt complexes $LW\Omega_{S}$, $W\Omega^{\ast}_{S}$, and $\mathcal{W}\Omega^{\ast}_{S}$ agree with each other. We refer the reader to \cite{ill} and \cite{blm} for detailed constructions of $W\Omega^{\ast}_{X}$ and $\mathcal{W}\Omega^{\ast}_{X}$ as well as relevant objects, e.g., Dieudonn\'e algebras \cite[Ch. 3]{blm}. \\
\indent Let us briefly recall the notion of (nonabelian) derived functors from \cite{htt}. Objects in the category $\CAlg^{\heartsuit}_{\bbF_{p}}$ of ordinary commutative $\bbF_{p}$-algebras can be described as filtered colimits of finitely presented $\bbF_{p}$-algebras, each of which is a (reflexive) coequalizer of finitely generated polynomial algebras. Reversing the process, one can consider the $\infty$-category $\CAlg^{\an}_{\bbF_{p}}$ of derived commutative $\bbF_{p}$-algebras which is freely generated by the category $\Poly_{\bbF_{p}}$ of finitely generated polynomial $\bbF_{p}$-algebras under sifted colimits; concretely, $\CAlg^{\an}_{\bbF_{p}}$ is equivalent to the full subcategory of space-valued presheaves on $\CAlg^{\poly}_{\bbF_{p}}$ spanned by finite product preserving functors. The category $\CAlg_{\bbF_{p}}^{\an}$ is equivalent to the underlying $\infty$-category of the simplicial model category of simplicial commutative $\bbF_{p}$-algebras \cite[Cor. 5.5.9.3]{htt}. By definition, for any $\infty$-category $\mathcal{D}$ admitting sifted colimits, composition with the Yoneda embedding $\Poly_{\bbF_{p}}\hookrightarrow\CAlg^{\an}_{\bbF_{p}}$ induces an equivalence $\Fun_{\Sigma}(\CAlg^{\an}_{\bbF_{p}},\mathcal{D})\simeq \Fun(\Poly_{\bbF_{p}},\mathcal{D})$ with inverse given by left Kan extension; here the source consists of sifted colimit preserving functors. Given a functor $F:\Poly_{\bbF_{p}}\to\mathcal{D}$, we call its left Kan extension $LF\in\Fun_{\Sigma}(\CAlg^{\an}_{\bbF_{p}},\mathcal{D})$ along the Yoneda embedding a (nonabelian) \emph{derived} functor of $F$. In particular, given a functor $G:\CAlg^{\heartsuit}_{\bbF_{p}}\to\mathcal{D}$, there is a natural map $L(G|_{\Poly_{\bbF_{p}}})|_{\CAlg^{\heartsuit}_{\bbF_{p}}}\to G$ in $\Fun(\CAlg^{\heartsuit}_{\bbF_{p}},\mathcal{D})$ which is up to equivalence a unique one extending the identity map of $G|_{\Poly_{\bbF_{p}}}$. 

\begin{lemma}\label{lem:derWr}
For each $r\geq 1$, there is an equivalence of functors $LW\Omega_{(-)}\dt_{\bbZ}\bbZ/p^{r}\simeq LW_{r}\Omega_{(-)}$ in $\Fun_{\Sigma}(\CAlg^{\an}_{\bbF_{p}},\mathcal{D}(\bbZ))$. 
\begin{proof}
Consider the functor $R\mapsto W\Omega_{R}^{\ast}\otimes_{\bbZ}\bbZ/p^{r}$ in $\Fun(\Poly_{\bbF_{p}},\mathcal{D}(\bbZ))$. For a smooth $\bbF_{p}$-algebra $R$ (e.g., $R\in\Poly_{\bbF_{p}}$), $W\Omega_{R}^{\ast}\otimes_{\bbZ}\bbZ/p^{r}\simeq W_{r}\Omega_{R}^{\ast}$ by \cite[Cor. I.3.17]{ill} and $p$-torsion freeness of each $W\Omega_{R}^{i}$ \cite[Cor. I.3.6]{ill}. On the other hand, the functor $R\mapsto LW\Omega_{R}\dt_{\bbZ}\bbZ/p^{r}$ in $\Fun(\CAlg^{\an}_{\bbF_{p}},\mathcal{D}(\bbZ))$ commutes with sifted colimits and for each $R\in\Poly_{\bbF_{p}}$, takes the value $W\Omega_{R}^{\ast}\dt_{\bbZ}\bbZ/p^{r}\simeq W\Omega_{R}^{\ast}\otimes_{\bbZ}\bbZ/p^{r}$ (using \cite[Cor. I.3.6]{ill}). Thus, both of the functors in question correspond to $W\Omega^{\ast}_{(-)}\otimes_{\bbZ}\bbZ/p^{r}$ (being its left Kan extensions), and hence are equivalent to each other. 
\end{proof} 
\end{lemma}

Now, let us specialize to the case of Cartier smooth rings. 

\begin{lemma}\label{lem:cartsmddR}
Let $S$ be a Cartier smooth $\bbF_{p}$-algebra. Then, the natural map $L\Omega_{S}\to\Omega^{\ast}_{S/\bbF_{p}}$ in $\mathcal{D}(\bbF_{p})$ is an equivalence. This equivalence respects conjugate filtrations, i.e., $\Fil^{\mathrm{conj}}_{\leq \sbullet}L\Omega_{S}\simeq \tau^{\leq \sbullet}\Omega^{\ast}_{S/\bbF_{p}}$. 
\begin{proof}
Here, $L\Omega_{S}$ is the derived de Rham complex of $S$, and the conjugate filtration $\Fil^{\mathrm{conj}}_{\leq\sbullet}L\Omega_{S}$ is defined as the value of the derived functor of the Postnikov filtration construction $R\mapsto \tau^{\leq \sbullet}\Omega^{\ast}_{R}$ at $S$. By construction, the natural map $\Fil^{\mathrm{conj}}_{\leq\sbullet}L\Omega_{S}\to\tau^{\leq\sbullet}\Omega^{\ast}_{S}$ enhances the natural map $L\Omega_{S}\to\Omega^{\ast}_{S}$ as a map of $\bbZ_{\geq 0}$-indexed increasing filtered objects, as Postnikov filtration and hence conjugate filtration admit $\Omega^{\ast}_{S}$ and $L\Omega_{S}$ respectively as their underlying objects. Thus, in order to show the claim, we need to check that the natural map $\Fil^{\mathrm{conj}}_{\leq i}L\Omega_{S}\to\tau^{\leq i}\Omega^{\ast}_{S}$ is an equivalence for each $i\geq 0$. By induction, it suffices to check that the induced natural map of $i$-th associated graded pieces is an equivalence for each $i\geq0$. By construction, $\gr_{i}\Fil^{\mathrm{conj}}_{\leq\sbullet}L\Omega_{S}\simeq L(R\mapsto H^{i}(\Omega^{\ast}_{R/\bbF_{p}})[-i])(S)\simeq L(R\mapsto \Omega^{i}_{R^{(1)}/\bbF_{p}}[-i])(S)\simeq L\Omega^{i}_{S^{(1)}}[-i]$ via Cartier isomorphism for polynomial algebras, and hence the natural map $\gr_{i}\Fil^{\mathrm{conj}}_{\leq\sbullet}L\Omega_{S}\to\gr_{i}(\tau^{\leq\sbullet}\Omega^{\ast}_{S/\bbF_{p}})$ between $i$-th associated graded pieces is given by the inverse Cartier map $L\Omega^{i}_{S^{(1)}}[-i]\simeq\left(\bigwedge^{i}L_{S^{(1)}/\bbF_{p}}\right)[-i]\to H^{i}(\Omega^{\ast}_{S/\bbF_{p}})[-i]$. The conditions required for Cartier smoothness of $S$ ensures this map is an equivalence for all $i\geq0$.  
\end{proof}
\end{lemma}

\begin{proposition}(\cite[Th. 9.4.1]{blm}) \label{prop:cartsmsaturated}
Let $S$ be a Cartier smooth $\bbF_{p}$-algebra. Then, the classical de Rham-Witt complex of \cite{ill} agrees with the saturated de Rham-Witt complex $\mathcal{W}\Omega^{\ast}_{S}$ of \cite{blm}. More precisely, the natural map $\gamma:W\Omega^{\ast}_{S}\to\mathcal{W}\Omega^{\ast}_{S}$ of commutative differential graded algebras in \cite[Cor. 4.4.11]{blm} is an isomorphism of strict Dieudonn\'e algebras. 
\begin{proof}
By \cite[Th. 9.4.1]{blm}, the map $\gamma_{1}:W_{1}\Omega^{\ast}_{S}=\Omega_{S/\bbF_{p}}^{\ast}\to\mathcal{W}_{1}\Omega^{\ast}_{S}$ is an isomorphism. Thus, it suffices to check the proof of \cite[Th. 4.4.12]{blm} which reduces the statement to the aforementioned proposition remains valid for $S$. It suffices to check $W\Omega^{\ast}_{S}$ is a saturated Dieudonn\'e algebra (which is in fact also strict). However, by $p$-torsion freeness \cite[Th. 2.8 (i)]{km} and $FW\Omega^{i}_{S} = d^{-1}(pW\Omega^{i+1}_{S})$ \cite[Th. 2.8 (ii)]{km} together with the injectivity of $F$, the claim follows. 
\end{proof}
\end{proposition}

By Proposition \ref{prop:cartsmsaturated} (which is in a sense a reinterpretation of \cite[Th. 2.8]{km}), we can freely use nice properties of $\mathcal{W}\Omega^{\ast}_{(-)}$ explained in \cite{blm} for $W\Omega^{\ast}_{S}$. For instance, the functor $R\mapsto \mathcal{W}\Omega^{\ast}_{R}:\CAlg_{\bbF_{p}}\to\text{DA}_{\text{str}}$ is a left adjoint functor \cite[Cor. 4.1.5]{blm}. Also, note that by \cite[Rem. 2.7.3]{blm} there are quasi-isomorphisms $W\Omega^{\ast}_{S}/p^{r} W\Omega^{\ast}_{S}\to W_{r}\Omega^{\ast}_{S}$ for each $r\geq 1$. We no longer distinguish between the two constructions for Cartier smooth algebras.

\begin{lemma}\label{lem:cartsmdWr}
Let $S$ be a Cartier smooth $\bbF_{p}$-algebra. Then, the natural map $LW_{r}\Omega_{S}\to W_{r}\Omega^{\ast}_{S}$ in $\mathcal{D}(\bbZ)$ is an equivalence for all $r\geq 1$. 
\begin{proof}
Let $R\in\Poly_{\bbF_{p}}$. By $p$-torsion freeness of each $W\Omega^{i}_{R}$ \cite[Cor. I.3.6]{ill}, we have an exact sequence of complexes $0\to W\Omega_{R}^{\ast}/p\xrightarrow{p^{r}}W\Omega_{R}^{\ast}/p^{r+1}\to W\Omega_{R}^{\ast}/p^{r}\to0$. This sequence induces a fiber sequence $W_{1}\Omega^{\ast}_{R}\to W_{r+1}\Omega^{\ast}_{R}\to W_{r}\Omega^{\ast}_{R}$ in $\mathcal{D}(\bbZ)$ by \cite[Cor. I.3.17]{ill}, which is natural in $R$. On the other hand, by Proposition \ref{prop:cartsmsaturated}, the argument remains valid for Cartier smooth algebras replacing $R$. Thus, there is a natural map of fiber sequences in $\mathcal{D}(\bbZ)$ from $LW_{1}\Omega_{S}\to LW_{r+1}\Omega_{S}\to LW_{r}\Omega_{S}$ to $W_{1}\Omega^{\ast}_{S}\to W_{r+1}\Omega^{\ast}_{S}\to W_{r}\Omega^{\ast}_{S}$. This reduces our problem to the case of $r=1$, i.e., to check the natural map $L\Omega_{S}\to \Omega^{\ast}_{S/\bbF_{p}}$ is an equivalence. This follows from Lemma \ref{lem:cartsmddR}. 
\end{proof}
\end{lemma}

Recall that each de Rham-Witt complex admits a $\bbZ_{\geq0}$-indexed descending filtration $\Nyg^{\geq\sbullet}W\Omega_{X}$ called \emph{Nygaard filtration}, which behaves as the $p$-adic filtration on the de Rham complex of the flat $p$-adic formal scheme with smooth special fiber. For Cartier smooth rings, the Nygaard filtration on $W\Omega^{\ast}_{S}$ takes the form $\Nyg^{\geq i}W\Omega^{\ast}_{S} = \cdots\xrightarrow{d}pVW\Omega^{i-2}_{S}\xrightarrow{d}VW\Omega^{i-1}_{S}\xrightarrow{d}W\Omega^{i}_{S}\xrightarrow{d}W\Omega^{i+1}_{S}\xrightarrow{d}\cdots$, where $V$ is the Verschiebung map of $W\Omega^{\ast}_{S}$. Along with the canonical inclusion map $\Nyg^{\geq i}W\Omega^{\ast}_{S}\to W\Omega^{\ast}_{S}$, recall that there is a natural map $\varphi/p^{i}:\Nyg^{\geq i}W\Omega^{\ast}_{S}\to W\Omega^{\ast}_{S}$ called the $i$-th \emph{divided Frobenius} map for each $i\geq0$ defined as follows. Write $F$ for the Frobenius map of $W\Omega^{\ast}_{S}$ as before; by construction it is an endomorphism of the underlying graded module interacting with the differential as $dF(x) = pF(dx)$. Thus, by taking $\varphi = p^{n}F:W\Omega^{n}_{S}\to W\Omega^{n}_{S}$ for each $n\in\bbZ$, we have an endomorphism $\varphi$ of the complex $W\Omega^{\ast}_{S}$. For each $i\geq 0$, the restriction of $\varphi$ on $\Nyg^{\geq i}W\Omega^{\ast}_{S}$ is divisible by $p^{i}$ and induces a well-defined map $\varphi/p^{i}:\Nyg^{\geq i}W\Omega^{\ast}_{S}\to W\Omega^{\ast}_{S}$ due to the $p$-torsion freeness of the target. 

\begin{lemma}(\cite[Lem. 8.2]{bms2} for smooth algebras) \label{lem:cartsmconj}
Let $S$ be a Cartier smooth $\bbF_{p}$-algebra. Then, the composition $\Nyg^{\geq i}W\Omega^{\ast}_{S}\xrightarrow{\varphi/p^{i}}W\Omega^{\ast}_{S}\xrightarrow{\text{mod}~p}\Omega^{\ast}_{S/\bbF_{p}}$ factors through $\tau^{\leq i}\Omega^{\ast}_{S/\bbF_{p}}$ and maps $\Nyg^{\geq i+1}W\Omega^{\ast}_{S}$ to zero. Moreover, the induced map $\gr^{i}_{\Nyg}W\Omega^{\ast}_{S} = \Nyg^{\geq i}W\Omega^{\ast}_{S}/\Nyg^{\geq i+1}W\Omega^{\ast}_{S}\to\tau^{\leq i}\Omega^{\ast}_{S/\bbF_{p}}$ is a quasi-isomorphism. 
\begin{proof}
In fact, proof \emph{loc. cit.} works due to \cite[Th. 2.8]{km}. Alternatively, by Proposition \ref{prop:cartsmsaturated}, $W\Omega^{\ast}_{S}$ is a saturated Dieudonn\'e complex, and hence \cite[Prop. 8.2.1]{blm} applies (where, notationally, $M^{\ast}\xrightarrow{\alpha_{F}}\eta_{p}M^{\ast}\to M^{\ast}$ equals the Frobenius $\varphi$ used here). Note that $\tau^{\leq i}(W\Omega^{\ast}_{S}/p)$ is quasi-isomorphic to $\tau^{\leq i}\Omega^{\ast}_{S/\bbF_{p}}$.
\end{proof}
\end{lemma}

\begin{proposition}\label{prop:cartsmLdW}
Let $S$ be a Cartier smooth $\bbF_{p}$-algebra. Then, the natural map $LW\Omega_{S}\to W\Omega_{S}^{\ast}$ in $\CAlg(\mathcal{D}(\bbZ))$ is an equivalence. This equivalence respects Nygaard filtrations and divided Frobenii. 
\begin{proof}
The functor $LW\Omega_{(-)}$ takes values in the $p$-complete derived category $\widehat{\mathcal{D}}(\bbZ_{p})$ in $\mathcal{D}(\bbZ)$, and hence a left adjoint $p$-completion together with Lemma \ref{lem:derWr} and Lemma \ref{lem:cartsmdWr} provide natural equivalences $LW\Omega_{S}\simeq \lim_{r}LW\Omega_{S}\dt_{\bbZ}\bbZ/p^{r}\simeq \lim_{r}LW_{r}\Omega_{S}\simeq \lim_{r}W_{r}\Omega^{\ast}_{S}$ in $\mathcal{D}(\bbZ)$. As each of the pro-objects $\{W_{r}\Omega^{i}_{S}\}_{r}$ has surjective transition maps, the last object is equivalent to $W\Omega^{\ast}_{S}$. On smooth $\bbF_{p}$-algebras, this construction reduces to the natural identification $LW\Omega_{R} = W\Omega^{\ast}_{R}$, and hence the composed natural equivalence is equivalent to the natural map $LW\Omega_{S}\to W\Omega^{\ast}_{S}$ of question. \\
\indent Consider the diagram of $\bbE_{\infty}$-rings in $\mathcal{D}(\bbZ)$ 
\begin{equation*}
\begin{tikzcd}
\Nyg^{\geq i}LW\Omega_{S} \arrow{r} \arrow{d}[swap]{\varphi_{i}} & \Nyg^{\geq i} W\Omega^{\ast}_{S} \arrow{d}{\varphi/p^{i}} \\ LW\Omega_{S} \arrow{r}[swap]{\simeq} & W\Omega^{\ast}_{S} 
\end{tikzcd}
\end{equation*}
induced from left Kan extensions (and similarly for the canonical inclusion maps in place of divided Frobenii $\varphi/p^{i}$). Here, we of course write $\Nyg^{\geq i}LW\Omega_{S}$ for $L(\Nyg^{\geq i}W\Omega_{(-)})(S)$. The top horizontal arrow for $i=0$ is the same as the bottom horizontal arrow by definition. To prove the top horizontal arrow is an equivalence for all $i\geq 0$, we consider the diagram
\begin{equation*}
\begin{tikzcd}
\Nyg^{\geq i+1}LW\Omega_{S} \arrow{r} \arrow{d} & \Nyg^{\geq i}LW\Omega_{S} \arrow{r} \arrow{d} & \Fil^{\text{conj}}_{\leq i} L\Omega_{S} \arrow{d}{\simeq}  \\   \Nyg^{\geq i+1} W\Omega^{\ast}_{S} \arrow{r} & \Nyg^{\geq i}W\Omega^{\ast}_{S} \arrow{r} & \tau^{\leq i}\Omega^{\ast}_{S/\bbF_{p}}
\end{tikzcd}
\end{equation*}
for each $i\geq 0$, where the horizontal lines are fiber sequences \cite[Lem. 8.2]{bms2}, Lemma \ref{lem:cartsmconj}. By Lemma \ref{lem:cartsmddR}, we know the right vertical arrow is an equivalence. Now, by induction on $i$ we have the claimed result. 
\end{proof}
\end{proposition}

\begin{remarkn}
Combined with Lemma \ref{lem:cartsmddR} and Lemma \ref{lem:cartsmconj}, Proposition \ref{prop:cartsmLdW} in particular gives an equivalence $\gr^{i}_{\Nyg}LW\Omega_{S}\simeq \Fil^{\conj}_{\leq i}L\Omega_{S}$ for Cartier smooth $S$. Note that in the language of prismatic complexes this takes the form $\gr^{i}_{\mathrm{Nyg}}F^{\ast}\prism_{S/\bbZ_{p}}\simeq\Fil^{\conj}_{\leq i}\overline{\prism}_{S}$; this equivalence itself does not require $S$ to be Cartier smooth, and the base (in this case, $\bbF_{p}$) being perfectoid is sufficient to guarantee this. 
\end{remarkn}

\indent Let $X$ be an $\bbF_{p}$-scheme, and view $W\Omega^{\ast}_{X}$ as a complex in $\Shv_{\Ab}(X_{\proet})$ (with an additional structure of a Dieudonn\'e complex). More precisely, for an $\bbF_{p}$-scheme $X$, each $W_{r}\Omega^{i}_{X}$ defined as the presheaf $(\Spec A)_{X}\mapsto W_{r}\Omega^{\ast}_{A}$ on $X_{\et}^{\aff}$ is a sheaf \cite[Prop. I.1.14]{ill} (cf. \cite[Th. 5.3.7]{blm}), and by taking the left adjoint geometric morphism $\nu^{\ast}:\Shv_{\Ab}(X_{\et})\to\Shv_{\Ab}(X_{\proet})$ we view (each $\nu^{\ast}W_{r}\Omega^{i}_{X}$ and) the limit $W\Omega^{i}_{X}$ naturally as a sheaf on $X_{\proet}$, cf. \cite[Prop. 5.6.2]{bs1}. As in \cite{bms2}, define $W_{r}\Omega_{X,\log}^{i}$ to be the image of the map $d\log:\bbG_{m/X}^{\otimes i}\to W_{r}\Omega^{i}_{X}=a_{1}\otimes\cdots\otimes a_{i}\mapsto d[a_{1}]/[a_{1}]\wedge\cdots\wedge d[a_{i}]/[a_{i}]$ of pro-\'etale sheaves on $X$, and let $W\Omega^{i}_{X,\log}$ be the (underived) limit of them. 

\begin{proposition}(\cite[Prop. 8.4]{bms2} for smooth schemes over perfect base) \label{prop:bms2log}
Let $X$ be a Cartier smooth $\bbF_{p}$-scheme. Then, the sequence of complexes of pro-\'etale sheaves $0\to W\Omega^{i}_{X,\log}[-i]\to \Nyg^{\geq i}W\Omega^{\ast}_{X}\xrightarrow{\varphi/p^{i}-1}W\Omega^{\ast}_{X}\to 0$ is exact, and $W\Omega^{i}_{X,\log} \simeq \lim_{r}W_{r}\Omega^{i}_{X,\log}$, where the limit is taken in $\mathcal{D}(\bbZ)$. 
\begin{proof}
Since each $W\Omega^{n}_{A}$ for affine open $\Spec A$ of $X$ is $p$-torsion free \cite[Th. 2.8 (i)]{km}, the divided Frobenius $\varphi/p^{i}:\Nyg^{\geq i}W\Omega^{\ast}_{X}\to W\Omega^{\ast}_{X}$ is defined as in the smooth case. Now, the proof in \cite[Prop. 8.4]{bms2} for smooth schemes over perfect fields works \emph{ad verbum}. Let us reproduce the proof \emph{loc. cit.} here for convenience. Let $\Spec A$ be an affine open of $X$. For degrees $n\neq i$, the map $(\Nyg^{\geq i}W\Omega^{\ast}_{A})^{i}\xrightarrow{\varphi/p^{i}-1}W\Omega^{i}_{A}$ is an isomorphism. Indeed, for $n>i$, $\varphi/p^{i}-1 = p^{n-i}F -1:W\Omega^{n}_{A}\to W\Omega^{n}_{A}$, and since $p^{n-i}F$ is $p$-adically contracting and $W\Omega^{n}_{A}$ is $p$-complete, this map is invertible. For $n<i$, consider the commutative diagram
\begin{equation*}
\begin{tikzcd}
W\Omega^{n}_{A} \arrow{r}{p^{i-1-n}V} \arrow[bend right=20, swap]{rr}{1-p^{i-1-n}V} & \Nyg^{\geq i}W\Omega^{n}_{A} \arrow{r}{\varphi/p^{i}-1} & W\Omega^{n}_{A}.  
\end{tikzcd}
\end{equation*}
The first map $p^{i-1-n}V$ is an isomorphism by definition. The lower curved arrow is also an isomorphism, since for $n<i-1$, $p^{i-1-n}V$ is $p$-adically contracting, and for $n=i-1$, $W\Omega^{n}_{A}$ is $V$-complete (which follows from \cite[Prop. 2.6 (ii)]{km}). Thus, the map $\varphi/p^{i}-1$ is an isomorphism. \\
\indent It remains to study the case of $n=i$, i.e., the map $W\Omega^{i}_{X}\xrightarrow{F-1}W\Omega^{i}_{X}$. By \cite[Cor. 4.1]{logHW} (extending \cite[Th. I.5.7.2]{ill} in the smooth case), the sequence $0\to\{W_{r}\Omega^{i}_{X,\log}\}_{r}\to\{W_{r}\Omega^{i}_{X}\}_{r}\xrightarrow{F-1}\{W_{r}\Omega^{i}_{X}\}_{r}\to0$ is exact in $\Pro(\Shv_{\Ab}(X_{\et}))$, and by taking the limit, we obtain the exact sequence $0\to W\Omega^{i}_{X,\log}\to W\Omega^{i}_{X}\xrightarrow{F-1}W\Omega^{i}_{X}\to 0$ in $\Shv_{\Ab}(X_{\proet})$ and $\lim^{q}_{r}W\Omega^{i}_{X,\log}\simeq 0$ for $q>0$. 
\end{proof}
\end{proposition}

\section{Syntomic cohomology of Cartier smooth rings}\label{sec:cartsmsyn}

Using the description of derived de Rham-Witt complexes of Cartier smooth rings from the previous section, we can immediately compute syntomic cohomology complexes $\bbZ_{p}(i)$ of Cartier smooth rings. Let us give a minimalistic exposition here sufficient for stating the results (especially in characteristic $p$), and refer the reader to \cite{bms2} for more details. Let $S$ be a $p$-quasisyntomic ring of characteristic $p$ (e.g., a Cartier smooth ring), and let $\CAlg^{\QSyn}_{S}$ be the full subcategory of the category of ordinary $S$-algebras consisting of $p$-quasisyntomic $S$-algebras $T$, i.e., whose structure map $S\to T$ is flat and $L_{T/S}$ has Tor amplitude in cohomological degrees $[-1,0]$. The opposite of this category admits a site structure with respect to the quasisyntomic topology via faithfully flat $p$-quasisyntomic maps, called the quasisyntomic site of $S$ \cite[Var. 4.35]{bms2}. Similarly, the (opposite of the) category $\CAlg^{\proet}_{S}$ of weakly \'etale $S$-algebras form the pro-\'etale site of $S$ \cite{bs1}. The derived de Rham-Witt complex construction $LW\Omega_{(-)}$ and its Nygaard filtration satisfy quasisyntomic descent on $\Spec S$, and hence the fiber $\bbZ_{p}(i)$ of $\varphi_{i}-\mathrm{can}:\Nyg^{\geq i}LW\Omega_{(-)}\to LW\Omega_{(-)}$ for each $i\geq0$ is a quasisyntomic sheaf on $\Spec S$. This fiber $\bbZ_{p}(i)$ precisely computes the $i$-th syntomic cohomology complex of $\Spec S$. 

\begin{corollary}(\cite[Cor. 8.21]{bms2} for smooth algebras over perfect fields)\label{cor:bms2log2}
Let $S$ be a Cartier smooth $\bbF_{p}$-algebra. Then, for each $i\geq 0$, the syntomic cohomology $\bbZ_{p}(i)$ in the pro-\'etale topology of $\Spec S$ is naturally equivalent to $W\Omega^{i}_{\Spec S,\log}[-i]$. More precisely, the pushforward of the sheaf $\bbZ_{p}(i)\in\Shv_{\Sp}((\CAlg_{S}^{\QSyn})^{\op})$ along the canonical map $(\CAlg^{\QSyn}_{S})^{\op}\to (\CAlg^{\proet}_{S})^{\op}$ of sites (induced from $\CAlg^{\proet}_{S}\subseteq\CAlg^{\QSyn}_{S}$) is naturally equivalent to $W\Omega^{i}_{\Spec S,\log}[-i]$. 
\begin{proof}
Note that $S$ is quasisyntomic. By definition $\bbZ_{p}(i)$ fits into the fiber sequence $\bbZ_{p}(i)\to \Nyg^{\geq i}LW\Omega_{-}\xrightarrow{\varphi_{i}-1}LW\Omega_{-}$ of sheaves on $\CAlg_{S}^{\QSyn}$, which after pushforward can be viewed as a fiber sequence in $\Shv_{\Sp}((\CAlg_{S}^{\proet})^{\op})\simeq\Shv_{\Sp}((\Spec S)_{\proet})$. By Lemma \ref{lem:frob} below, objects of $\CAlg_{S}^{\proet}$ are Cartier smooth. Hence, the map $\Nyg^{\geq i}LW\Omega_{-}\xrightarrow{\varphi_{i}-1}LW\Omega_{-}$ is equivalent to $\Nyg^{\geq i}W\Omega^{\ast}_{\Spec S}\xrightarrow{\varphi/p^{i}-1}W\Omega^{\ast}_{\Spec S}$ by Proposition \ref{prop:cartsmLdW}. Now, Proposition \ref{prop:bms2log} implies that there is a natural equivalence $\bbZ_{p}(i)_{(\Spec S)_{\proet}}\simeq W\Omega^{i}_{\Spec S,\log}[-i]$. 
\end{proof}
\end{corollary}

\begin{lemma}\label{lem:frob}
Let $R\to R'$ be a weakly \'etale map of $\bbF_{p}$-algebras. Following the convention of \cite[Not. 9.5.8]{blm}, set $R^{(1)} = R$ and write $R$ when we view it as an $R^{(1)}$-algebra via Frobenius $R^{(1)}\to R$. Then, \\
(1) The commutative diagram 
\begin{equation*}
\begin{tikzcd}
R^{(1)} \arrow{r}{F_{R}} \arrow{d} & R \arrow{d} \\
R'^{(1)} \arrow[swap]{r}{F_{R'}} & R'
\end{tikzcd}
\end{equation*}
is a pushout square of $\bbF_{p}$-algebras. \\
(2) If $R$ is Cartier smooth, then so is $R'$. 
\begin{proof}
(1) The relative Frobenius map $\Spec R'\to\Spec R'^{(1)}\times_{\Spec R^{(1)}}\Spec R$ is weakly \'etale by assumption \cite[Prop. 2.3.3 (4)]{bs1}. On the other hand, it is a universal homeomorphism \cite[tag 0CCB]{stacks}, and hence must be an isomorphism \cite[tag 0F6V]{stacks}.\\
\noindent (2) From $L_{R'^{(1)}/R^{(1)}}\simeq 0$, one observes $L_{R'^{(1)}/\bbF_{p}}$ is equivalent to a flat $R'^{(1)}$-module. As in the proof of \cite[Prop. 9.5.11]{blm}, (1) implies that the inverse Cartier operator $\Omega^{\sbullet}_{R'^{(1)}}\to H^{\sbullet}(\Omega^{\ast}_{R'})$ is an isomorphism of graded $R'^{(1)}$-algebras by base change.  
\end{proof}
\end{lemma}

\begin{corollary}\label{cor:bms2log3}
Let $S$ be a Cartier smooth $\bbF_{p}$-algebra. Then, for each $i\geq 0$ and $r\geq 1$, the mod-$p^{r}$ syntomic cohomology $\bbZ/p^{r}(i)$ in the \'etale topology of $\Spec S$ is naturally equivalent to $W_{r}\Omega^{i}_{\Spec S,\log}[-i]$. 
\begin{proof}
The equivalence is obtained by smashing the equivalence $\bbZ_{p}(i)\simeq W\Omega^{i}_{\Spec S,\log}[-i]$ of Corollary \ref{cor:bms2log2} (pushed forward to $\Shv_{\Sp}((\CAlg_{S}^{\et})^{\op})$) with the mod-$p^{r}$ Moore spectrum $\bbS/p^{r}$ pointwisely. By definition $\bbZ_{p}(i)/p^{r} = \bbZ/p^{r}(i)$. On the other hand, \cite[Th. 2.11]{km} implies that $W\Omega^{i}_{R,\log}/p^{r}\simeq W_{r}\Omega^{i}_{R,\log}$ naturally for any local Cartier smooth $\bbF_{p}$-algebra $R$. As Cartier smoothness is closed under filtered colimits, we know $W\Omega^{i}_{\Spec S,\log}/p^{r}\simeq W_{r}\Omega^{i}_{\Spec S,\log}$ from \'etale-local stalk computations. 
\end{proof}
\end{corollary}

\section{$\TC(-,\bbZ_{p})$ of Cartier smooth rings}\label{sec:cartsmTC}

Recall that for any ring $R$ of characteristic $p$, its topological Hochschild homology $\THH(R)$ is a connective and $p$-complete cyclotomic spectrum, and hence $\TC(R) = \TC(\THH(R))$ is $p$-complete. Let $S$ be a $p$-quasisyntomic ring of characteristic $p$. We know that $p$-adic $\TC$, i.e., the functor $\TC(-,\bbZ_{p}) = R\mapsto \TC(R)^{\wedge_{p}} = \TC(R,\bbZ_{p})$, satisfies $p$-complete faithfully flat descent \cite[Cor. 3.4]{bms2}, and in particular have $\TC(-)\simeq \TC(-,\bbZ_{p})\in\Shv_{\Sp}((\CAlg^{\QSyn}_{S})^{\op})$ and hence $\TC/p^{r}\in\Shv_{\Sp}((\CAlg^{\QSyn}_{S})^{\op})$ for each $r\geq 1$ as well. Through the construction of \cite[Sec. 7.4]{bms2}, we know that $\TC(-,\bbZ_{p})$ admits a $\bbZ_{\geq 0}$-indexed descending complete and exhaustive filtration, called the \emph{motivic} filtration, whose $i$-th associated graded piece is given by the (shift of the) $i$-th syntomic complex $\bbZ_{p}(i)_{\Spec S}[2i]$. Of course, the construction fundamentally relies on the quasisyntomic descent and the novel interpretation of $p$-adic $\TC$ by \cite{ns} in terms of the equalizer of the Frobenius and the canonical map between $p$-adic $\TC^{-}$ and $\TP$.  \\
\indent For a quasisyntomic $\bbF_{p}$-algebra $S$, let $\Shv_{\Sp}((\CAlg^{\QSyn}_{S})^{\op})\xrightarrow{\lambda_{\ast}}\Shv_{\Sp}((\CAlg_{S}^{\proet})^{\op})$ be the natural map induced from the geometric morphism comparing underlying $\infty$-topoi. Similarly, let $\Shv_{\Sp}((\CAlg^{\QSyn}_{S})^{\op})\xrightarrow{\lambda'_{\ast}}\Shv_{\Sp}((\CAlg_{S}^{\et})^{\op})$ be the analogous natural map. 

\begin{corollary}\label{cor:cartsmTCpadic}
Let $S$ be a Cartier smooth $\bbF_{p}$-algebra. Then, there is an isomorphism 
\begin{equation*}
\pi_{i}\lambda_{\ast}\TC(-,\bbZ_{p})\isomto W\Omega^{i}_{\Spec S,\log}
\end{equation*} 
in $\Shv_{\Ab}((\CAlg_{S}^{\proet})^{\op})$ for all $i\geq 0$, and the sheaf $\lambda_{\ast}\TC(-,\bbZ_{p})$ is connective. Here, $\pi_{i} = \pi_{i}^{\proet}$ is taken with respect to the canonical t-structure on $\Shv_{\Sp}((\CAlg_{S}^{\proet})^{\op})$ \cite[1.3.2.3]{sag}. 
\begin{proof}
By Proposition \ref{prop:filtspsq} applied to the image of the motivic filtration on $\TC(-,\bbZ_{p})$ by $\lambda_{\ast}$, there is a spectral sequence $E^{i,j}_{2} = \pi_{-i-j}(\lambda_{\ast}\bbZ_{p}(-j)[-2j])$ in $\Shv_{\Ab}((\CAlg_{S}^{\proet})^{\op})$ which converges (conditionally) to $\pi_{-i-j}\lambda_{\ast}\TC(-,\bbZ_{p})$. Corollary \ref{cor:bms2log2} tells us that the spectral sequence degenerates at the second page, and that $\pi_{i}\lambda_{\ast}\TC(-,\bbZ_{p})\simeq H^{i}(\lambda_{\ast}\bbZ_{p}(i))$, which is isomorphic to $W\Omega^{i}_{\Spec S,\log}$ for $i\geq0$ and is zero for $i<0$. 
\end{proof}
\end{corollary}

\begin{corollary}\label{cor:cartsmTCmodpr}
Let $S$ be a Cartier smooth $\bbF_{p}$-algebra. Then, there is an isomorphism 
\begin{equation*}
\pi_{i}\lambda'_{\ast}\TC/p^{r}\isomto W_{r}\Omega^{i}_{\Spec S,\log}
\end{equation*} 
in $\Shv_{\Ab}((\CAlg_{S}^{\et})^{\op})$ for all $i\geq 0$ and $r\geq 1$, and the sheaf $\lambda'_{\ast}\TC/p^{r}$ is connective. Here, $\pi_{i} = \pi_{i}^{\et}$ is taken with respect to the canonical t-structure on $\Shv_{\Sp}((\CAlg_{S}^{\et})^{\op})$ \cite[1.3.2.3]{sag}. 
\begin{proof}
From the motivic filtration $\cdots\to \Fil^{\geq 1}\TC(-,\bbZ_{p})\to\Fil^{\geq 0}\TC(-,\bbZ_{p}) = \TC(-,\bbZ_{p})$, we get the exhaustive complete $\bbZ_{\geq 0}$-indexed filtration $\cdots\to (\Fil^{\geq 1}\TC(-,\bbZ_{p}))/p^{r}\to(\Fil^{\geq 0}\TC(-,\bbZ_{p}))/p^{r}\simeq \TC/p^{r}$ on $\TC/p^{r}$, i.e., we set $\Fil^{\geq i}\TC/p^{r} = (\Fil^{\geq i}\TC(-,\bbZ_{p}))/p^{r}$. In particular, its associated graded object is $\gr^{i}\TC/p^{r} = \bbZ/p^{r}(i)$. Applying the right adjoint functor $\lambda'_{\ast}$, we obtain the exhaustive complete $\bbZ_{\geq 0}$-indexed filtration $\lambda'_{\ast}\Fil^{\geq i}\TC/p^{r}$. By applying Proposition \ref{prop:filtspsq} to this filtration, we have a spectral sequence $E^{i,j}_{2} = \pi_{-i-j}(\lambda'_{\ast}\bbZ/p^{r}(-j)[-2j])$ which (conditionally) converges to $\pi_{-i-j}\lambda'_{\ast}\TC/p^{r}$ in $\Shv_{\Ab}((\CAlg_{S}^{\et})^{\op})$. Corollary \ref{cor:bms2log3} tells us that the spectral sequence degenerates at the second page, and that $\pi_{i}\lambda'_{\ast}\TC/p^{r}\simeq H^{i}(\lambda'_{\ast}\bbZ/p^{r}(i))$, which is isomorphic to $W_{r}\Omega^{i}_{\Spec S,\log}$ for $i\geq0$ and is zero for $i<0$. 
\end{proof}
\end{corollary}

From now on, let us drop $\lambda_{\ast}$ (resp. $\lambda'_{\ast}$) and view $\TC(-,\bbZ_{p})$ (resp. $\TC/p^{r}$) as a pro-\'etale sheaf (resp. an \'etale sheaf) of spectra on $\Spec S$.  

\begin{proposition}\label{prop:cartsmKmodpr}
Let $S$ be a Cartier smooth $\bbF_{p}$-algebra. Also, let $\pi_{i}^{\et}(\K_{\geq 0}/p^{r})$ be the \'etale sheafification of the presheaf $\pi_{i}(\K_{\geq 0}(-)/p^{r})$ on $(\Spec S)_{\et}$. Then, there is an isomorphism 
\begin{equation*}
\pi^{\et}_{i}(\K_{\geq 0}/p^{r})\isomto W_{r}\Omega^{i}_{\Spec S,\log}
\end{equation*} 
for each $i\geq 0$ and $r\geq 1$. 
\begin{proof}
Consider the composition $\pi^{\et}_{i}(\K_{\geq 0}/p^{r})\xrightarrow{\tr}\pi_{i}(\TC/p^{r})\to W_{r}\Omega^{i}_{\Spec S,\log}$. The second map is the isomorphism of Corollary \ref{cor:cartsmTCmodpr}. To prove that the first map is also an isomorphism, it suffices to check the map induces isomorphisms on stalks. As $\K_{\geq 0}$ and $\TC/p^{r}$ commutes with filtered colimits of rings (see \cite[Th. 1.12]{cmm} for the latter), we are reduced to the statement that the map $\pi_{i}(\K_{\geq 0}(R)/p^{r})\xrightarrow{\tr}\pi_{i}(\TC(R)/p^{r})$ is an isomorphism for $R$ strictly Henselian local of characteristic $p$, and this holds by \cite[Th. 1.7]{cmm}. 
\end{proof}
\end{proposition}

In the next section, we will refine this to Zariski topology in the form of Corollary \ref{cor:cartsmglmodpr}. For the case of smooth rings over characteristic $p$ perfect fields, this was proved by Geisser-Levine \cite{gl}. 

\section{Geisser-Levine theorem for Cartier smooth rings via syntomic cohomology} \label{sec:cartsmglmotivic}

In this section, we complete our proof of Theorem \ref{th:cartsmglpadic} and deduce Corollary \ref{cor:cartsmglmodpr} as stated in \cite[Th. 2.1]{km}. As noted in Remark \ref{rmk:padicmodprequiv}, these two are equivalent statements. We also remark a quick consequence of the theorem to the higher algebraic K-groups of perfect rings. \\
\indent We first compute the descent spectral sequence of $\TC(-,\bbZ_{p})$ in pro-\'etale topology in Proposition \ref{prop:proetpostnikov}. Note that in pro-\'etale $\infty$-topos of schemes, hypercomplete objects are Postnikov complete \cite[Prop. 3.2.3]{bs1}. Nevertheless, we provide a different proof verifying the Postnikov completeness of $\TC(-,\bbZ_{p})$ and computing the resulting descent spectral sequence.   

\begin{lemma}\label{lem:indetdWr}
Let $f:\Spec S'\to\Spec S$ be the map given by an ind-\'etale map $S\to S'$ of commutative $\bbF_{p}$-algebras. Then, for each $r\geq 1$ and $j\geq 0$, there is a natural isomorphism $f^{\ast}\mathcal{W}_{r}\Omega^{\ast}_{\Spec S}\simeq \mathcal{W}_{r}\Omega^{\ast}_{\Spec S'}$ of \'etale sheaves of commutative differential graded algebras. 
\begin{proof}
Consider the square 
\begin{equation*}
\begin{tikzcd}
\Shv_{\Ab}((\Spec W_{r}(S))_{\et}) \arrow{r} \arrow{d}[swap]{f_{r}^{\ast}}& \Shv_{\Ab}((\Spec S)_{\et})\arrow{d}{f^{\ast}}\\
\Shv_{\Ab}((\Spec W_{r}(S'))_{\et}) \arrow{r}& \Shv_{\Ab}((\Spec S')_{\et})
\end{tikzcd}
\end{equation*}
induced from the commutative square
\begin{equation*}
\begin{tikzcd}
W_{r}(S) \arrow{r} \arrow{d}& S \arrow{d}\\
W_{r}(S') \arrow{r}& S'.
\end{tikzcd}
\end{equation*}
By topological invariance of the \'etale site, the horizontal arrows induced by pullbacks are equivalences. Moreover, under this equivalence, each $\mathcal{W}_{r}\Omega^{j}_{\Spec S}$ corresponds to the quasicoherent sheaf on $(\Spec W_{r}(S))_{\et}$ given by the $W_{r}(S)$-module $\mathcal{W}_{r}\Omega^{j}_{S}$, and similarly for $\mathcal{W}_{r}\Omega^{j}_{\Spec S'}$ \cite[Th. 5.3.7]{blm}. Thus, it suffices to check that the natural map $W_{r}(S')\otimes_{W_{r}(S)}\mathcal{W}_{r}\Omega^{\ast}_{S}\to \mathcal{W}_{r}\Omega^{\ast}_{S'}$ is an isomorphism of graded $W_{r}(S')$-algebras. Since $W_{r}(-)$ commutes with filtered colimits, by writing $S\to S'\simeq \colim_{\alpha}S'_{\alpha}$ as a filtered colimit of \'etale $S$-algebras, we have $W_{r}(S')\simeq\colim_{\alpha}W_{r}(S'_{\alpha})$ compatible with maps to $S'\simeq\colim_{\alpha}S'_{\alpha}$. Hence, $W_{r}(S')\otimes_{W_{r}(S)}\mathcal{W}_{r}\Omega^{\ast}_{S}\simeq \colim_{\alpha}W_{r}(S'_{\alpha})\otimes_{W_{r}(S)}\mathcal{W}_{r}\Omega^{\ast}_{S}\simeq \colim_{\alpha}\mathcal{W}_{r}\Omega^{\ast}_{S'_{\alpha}}\simeq \mathcal{W}_{r}\Omega^{\ast}_{S'}$, where the second isomorphism uses \cite[Cor. 5.3.5]{blm}, while the final isomorphism follows from \cite[Cor. 4.3.5]{blm}.
\end{proof}
\end{lemma}

\begin{proposition}\label{prop:proetpostnikov}
Let $S$ be a Cartier smooth $\bbF_{p}$-algebra. Then, the followings hold:\\
(1) The prestable $\infty$-category of connective pro-\'etale sheaves $\Shv_{\Sp}((\Spec S)_{\proet})_{\geq 0}$ has enough objects of cohomological dimension $\leq 1$ with coefficients in $\{\pi_{j}\TC(-,\bbZ_{p})\}_{j\in\bbZ}$ (cf. \cite[Def. 2.8]{etalek}). \\
(2) The pro-\'etale sheaf $\TC(-,\bbZ_{p})$ is Postnikov complete, and the resulting descent spectral sequence gives functorial exact sequences $0\to H^{1}((\Spec S)_{\proet}, W\Omega^{i+1}_{\Spec S,\log})\to \pi_{i}\TC(S,\bbZ_{p})\to W\Omega^{i}_{S,\log}\to0$ for each $i\geq0$.
\begin{remark}
Similarly, for each $r\geq 1$, the \'etale descent spectral sequence for $\pi_{j}(\TC/p^{r})$ gives functorial exact sequences $0\to H^{1}((\Spec S)_{\et},W_{r}\Omega^{i+1}_{\Spec S,\log})\to\pi_{i}(\TC(S)/p^{r})\to W_{r}\Omega^{i}_{S,\log}\to0$ for each $i\geq 0$ using Corollary \ref{cor:cartsmTCmodpr} (cf. \cite[Cor. 5.20]{etalek}).  
\end{remark}
\begin{proof}
(1) Note that the category is equivalent to $\Shv_{\Sp^{\cn}}((\Spec S)_{\proet})$ \cite[1.3.5.7]{sag}, and hence is generated under colimits by representable sheaves. It suffices to check that for each $j\geq 0$ and weakly \'etale map $\Spec S'\to \Spec S$, there exists a faithfully flat weakly \'etale map $\Spec S''\to\Spec S'$ such that $H^{i}((\Spec S'')_{\proet},W\Omega^{j}_{\Spec S,\log})=0$ for $i>1$ by Corollary \ref{cor:cartsmTCpadic}. Since weakly \'etale maps can be refined by ind-\'etale maps \cite[Th. 2.3.4]{bs1}, we are reduced to check that for each $j\geq 0$ and a map $\Spec S'\to \Spec S$ given by an ind-\'etale ring map $S\to S'$, we have $H^{i}((\Spec S')_{\proet}, W\Omega^{j}_{\Spec S,\log})=0$ for $i>1$. The exact sequence $0\to W\Omega^{j}_{\Spec S,\log}\to W\Omega^{j}_{\Spec S}\xrightarrow{F-1}W\Omega^{j}_{\Spec S}\to0$ of Proposition \ref{prop:bms2log} further reduces this to the vanishing $H^{i}((\Spec S')_{\proet},W\Omega^{j}_{\Spec S}|_{\Spec S'})=0$ for $i>0$. Note that $W\Omega^{j}_{\Spec S}|_{\Spec S'}\simeq (\Rlim_{r}\nu^{\ast}W_{r}\Omega^{j}_{\Spec S})|_{\Spec S'}\simeq \Rlim_{r}\nu^{\ast}(W_{r}\Omega^{j}_{\Spec S}|_{\Spec S'})$ using \cite[Lem. 5.4.1]{bs1} and that $\Spec S'\to\Spec S$ is weakly \'etale. We then have
\begin{align*} 
&\RG((\Spec S')_{\proet},W\Omega^{j}_{\Spec S}|_{\Spec S'})\simeq \RG((\Spec S')_{\proet},\Rlim_{r}\nu^{\ast}(W_{r}\Omega^{j}_{\Spec S}|_{\Spec S'}))\\
&\simeq \Rlim_{r}\RG((\Spec S')_{\proet}, \nu^{\ast}(W_{r}\Omega^{j}_{\Spec S}|_{\Spec S'}))\simeq \Rlim_{r}\RG((\Spec S')_{\et}, W_{r}\Omega^{j}_{\Spec S}|_{\Spec S'})\\
&\simeq \Rlim_{r}\RG((\Spec S')_{\et}, W_{r}\Omega^{j}_{\Spec S'})\simeq \Rlim_{r}W_{r}\Omega^{j}_{S'}\simeq W\Omega^{j}_{S'},
\end{align*}
where the second equivalence uses the fact that $\RG$ and $\Rlim$ commutes, the third equivalence uses \cite[Cor. 5.1.6]{bs1}, the 4th equivalence follows from Lemma \ref{lem:indetdWr}, and the 5th equivalence follows from \cite[Th. 5.3.7]{blm} and its proof. This in particular implies $H^{0}((\Spec S')_{\proet}, W\Omega^{j}_{\Spec S}|_{\Spec S'}) = W\Omega^{j}_{S'}$ and the claimed vanishing for nonzero degrees. \\
(2) Recall that the quasisyntomic sheaf $\TC(-,\bbZ_{p})$ is hypercomplete. Since the left adjoint geometric morphism $\Shv_{\Sp}((\Spec S)_{\proet})\to\Shv_{\Sp}((\CAlg^{\QSyn}_{S})^{\op})$ is t-exact (in particular preserves $\infty$-connective objects), the corresponding right adjoint geometric morphism (i.e., the pushforward) preserves hypercompleteness, and in particular $\TC(-,\bbZ_{p})$ understood as a pro-\'etale sheaf is hypercomplete. Then, by (1) above and \cite[Prop. 2.10]{etalek}, we know $\TC(-,\bbZ_{p})$ is in fact a Postnikov complete pro-\'etale sheaf. In particular, we have a descent spectral sequence $E^{i,j}_{2}=H^{i}((\Spec S)_{\proet},W\Omega^{j}_{\Spec S,\log})$ which converges conditionally to $\pi_{-i+j}\TC(S,\bbZ_{p})$. By 
\begin{equation*}
\RG((\Spec S)_{\proet},W\Omega^{j}_{\Spec S})\simeq W\Omega^{j}_{S}
\end{equation*} 
computed above in the process of proving (1), we know $E_{2}^{i,j}=0$ for $i\neq 0,1$. Thus, the spectral sequence degenerates at the second page, and we obtain the claimed exact sequence.  
\end{proof} 
\end{proposition}

Now, we are ready to prove the main theorem of this section:

\begin{theorem}\label{th:cartsmglpadic}
Let $S$ be a local Cartier smooth $\bbF_{p}$-algebra. Then, the map $\pi_{i}\K_{\geq 0}(S,\bbZ_{p})\xrightarrow{\tr}\pi_{i}\TC(S,\bbZ_{p})\to W\Omega^{i}_{S,\log}$ given as a composition of the cyclotomic trace and the natural map of Propostion \ref{prop:proetpostnikov} is an isomorphism for all $i\geq 0$. 
\begin{proof}
Consider the map 
\begin{equation*}
\map_{\Shv_{\Sp}((\Spec S)_{\proet})}(h_{(\Spec S)_{\proet}},\tau_{\leq\ast}\TC(-,\bbZ_{p}))\to\map_{\Shv_{\Sp}((\Spec S)_{\et})}(h_{(\Spec S)_{\et}},\tau_{\leq\ast}(\TC/p^{r}))
\end{equation*}
in $\Fun(\bbZ,\Sp)$ induced from the right adjoint geometric morphism $\nu_{\ast}$(together with the canonical map $h_{(\Spec S)_{\et}}\to\nu_{\ast}h_{(\Spec S)_{\proet}}$ and with taking mod-$p^{r}$ cofiber). This induces a map between the associated descent spectral sequences, and in particular by Proposition \ref{prop:proetpostnikov} induces a commutative square
\begin{equation*}
\begin{tikzcd}
\pi_{i}(\TC(S,\bbZ_{p})) \arrow{r} \arrow{d} & \pi_{i}(\TC(S)/p^{r})\arrow{d} \\
W\Omega^{i}_{S,\log} \arrow{r}[swap]{\text{can}} & W_{r}\Omega^{i}_{S,\log}.
\end{tikzcd}
\end{equation*} 
Thus, the composition $\pi_{i}\TC(S,\bbZ_{p})\to \lim_{r}\pi_{i}(\TC(S)/p^{r})\to \lim_{r}W_{r}\Omega^{i}_{S,\log}$ agrees with the map $\pi_{i}\TC(S,\bbZ_{p})\to W\Omega^{i}_{S,\log}$ from the pro-\'etale descent spectral sequence. Taking cyclotomic trace maps into account, we have a commutative diagram
\begin{equation*}
\begin{tikzcd}
\pi_{i}(\K_{\geq 0}(S,\bbZ_{p})) \arrow{r} \arrow{d}[swap]{\tr} & \lim_{r}\pi_{i}(\K_{\geq 0}(S)/p^{r}) \arrow{d}{\tr} \arrow[bend left=80]{dd} \\
\pi_{i}(\TC(S,\bbZ_{p})) \arrow{r} \arrow{rd} & \lim_{r}\pi_{i}(\TC(S)/p^{r}) \arrow{d}\\
 & W\Omega^{i}_{S,\log}.
\end{tikzcd}
\end{equation*}
Let $R$ be a local ind-smooth $\bbF_{p}$-algebra, and set $S=R$ in the diagram above. By Geisser-Levine \cite[Th. 8.3]{gl}, we know each $\pi_{i}(\K_{\geq 0}(R))$ is $p$-torsion free and that the cyclotomic trace induces isomorphisms $\pi_{i}(\K_{\geq 0}(R)/p^{r})\xrightarrow{\tr}\pi_{i}(\TC(R)/p^{r})\to W_{r}\Omega^{i}_{R,\log}$ compatible with each other for all $r\geq 1$. Here, the second map is obtained from the \'etale descent spectral sequence for $\TC/p^{r}$. Hence, the top horizontal map (due to $p$-torsion freeness) and the rightmost bent arrow are isomorphisms, i.e., we have the claimed result for $S=R$ ind-smooth over $\bbF_{p}$. In particular, $\pi_{i}(\K_{\geq 0}(R,\bbZ_{p}))\xrightarrow{\tr}\pi_{i}(\TC(R,\bbZ_{p}))$ is split injective.\\
\indent Now, we consider $S$ as stated above. Let $R$ be an ind-smooth local $\bbF_{p}$-algebra which admits a surjection onto $S$ with Henselian kernel \cite{hoo}. Consider the commutative diagram
\begin{equation*}
\begin{tikzcd}
\pi_{i}(\K_{\geq 0}(R,\bbZ_{p})) \arrow{r}{\tr} \arrow{d} & \pi_{i}(\TC(R,\bbZ_{p})) \arrow{r} \arrow{d} & W\Omega^{i}_{R,\log} \arrow{d} \\
\pi_{i}(\K_{\geq 0}(S,\bbZ_{p})) \arrow{r}{\tr} & \pi_{i}(\TC(S,\bbZ_{p})) \arrow{r} & W\Omega^{i}_{S,\log}. 
\end{tikzcd}
\end{equation*}
Here, the right horizontal maps are obtained by Proposition \ref{prop:proetpostnikov}. By rigidity \cite[Th. 1.5]{cmm} and the (split) injectivity of $\tr$ for $R$, we have short exact Mayer-Vietoris type sequences $0\to\pi_{i}(\K_{\geq 0}(R,\bbZ_{p}))\to\pi_{i}(\TC(R,\bbZ_{p}))\oplus\pi_{i}(\K_{\geq 0}(S,\bbZ_{p}))\to\pi_{i}(\TC(S,\bbZ_{p}))\to 0$, i.e., the left square in the diagram is biCartesian (in particular a pushout square). To check that the right square in the diagram is biCartesian, it suffices to check the map between kernels of horizontal maps of the right square is an isomorphism. By Proposition \ref{prop:proetpostnikov}, this map equals $H^{1}((\Spec R)_{\proet},W\Omega^{i+1}_{\Spec R,\log})\to H^{1}((\Spec S)_{\proet},W\Omega^{i+1}_{\Spec S,\log})$, which can be rewritten as the natural map $H^{i+2}(\bbZ_{p}(i+1)(R))\to H^{i+2}(\bbZ_{p}(i+1)(S))$ using Corollary \ref{cor:cartsmTCpadic}. By rigidity \cite[Th. 5.2]{ammn} of syntomic cohomology, we know this map is an isomorphism.  
\end{proof}
\end{theorem}

\begin{corollary}(\cite[Th. 2.1]{km})\label{cor:cartsmglmodpr}
Let $S$ be a local Cartier smooth $\bbF_{p}$-algebra. Then, each $\pi_{i}\K_{\geq0}(S)$ is $p$-torsion free and the map $\pi_{i}(\K_{\geq 0}(S)/p^{r})\xrightarrow{\tr}\pi_{i}(\TC(S)/p^{r})\to W_{r}\Omega^{i}_{S,\log}$ induced by the cyclotomic trace and the natural map is an isomorphism for all $i\geq 0$ and $r\geq 1$. 
\begin{proof}
Consider the following diagram
\begin{equation*}
\begin{tikzcd}
& \K^{M}_{i}(S)/p^{r} \arrow{d} \arrow[bend right=80, two heads, swap]{ddd}{d\log} &&&\\
0\arrow{r} & \pi_{i}(\K_{\geq 0}(S))/p^{r}\arrow{r} \arrow[hook]{d} & \pi_{i}(\K_{\geq 0}(S)/p^{r})\arrow{r} \arrow{d}{\simeq} & \pi_{i-1}(\K_{\geq 0}(S))[p^{r}]\arrow{r} \arrow[two heads]{d}&0 \\
0\arrow{r} & \pi_{i}(\K_{\geq 0}(S,\bbZ_{p}))/p^{r}\arrow{r} \arrow{d}{\simeq} & \pi_{i}(\K_{\geq 0}(S,\bbZ_{p})/p^{r})\arrow{r} & \pi_{i-1}(\K_{\geq 0}(S,\bbZ_{p}))[p^{r}]\arrow{r} &0.\\
& W_{r}\Omega^{i}_{S,\log} &&&
\end{tikzcd}
\end{equation*}
The two exact rows are derived from exact sequences induced from mod $p^{r}$ fiber sequences for $\K_{\geq0}(S)$ and $\K_{\geq 0}(S,\bbZ_{p})$ respectively, and in particular the two squares in the diagram commute. By Theorem \ref{th:cartsmglpadic} and Proposition \ref{prop:cartsmsaturated}, the bottom right object $\pi_{i-1}(\K_{\geq 0}(S,\bbZ_{p}))[p^{r}]$ is zero. By the usual diagram chasing, we know that among those squares the left vertical arrow is injective and has cokernel isomorphic to $\pi_{i-1}(\K_{\geq 0}(S))[p^{r}]$. \\
\indent The left bottom vertical isomorphism is the map induced by the isomorphism of Theorem \ref{th:cartsmglpadic}, and hence by the commutative diagram
\begin{equation*}
\begin{tikzcd}
\pi_{i}(\K_{\geq 0}(S))/p^{r} \arrow{r} \arrow[swap]{d}{\tr} & \pi_{i}(\K_{\geq 0}(S,\bbZ_{p}))/p^{r}\arrow{d}{\tr}\\
\pi_{i}(\TC(S))/p^{r} \arrow{r} \arrow{d} & \pi_{i}(\TC(S,\bbZ_{p}))/p^{r}\arrow{d} \arrow{ld}\\
\pi_{i}(\TC(S)/p^{r}) \arrow{r} & W_{r}\Omega^{i}_{S,\log}
\end{tikzcd}
\end{equation*}
the map $\pi_{i}(\K_{\geq 0}(S))/p^{r}\to\pi_{i}(\K_{\geq 0}(S,\bbZ_{p}))/p^{r}\to W_{r}\Omega^{i}_{S,\log}$ in the first diagram equals the map obtained as a composition of the cyclotomic trace and the natural map obtained from the \'etale descent spectral sequence for $\TC/p^{r}$. Since cyclotomic trace lifts the $d\log$ map on symbols \cite[Lem. 4.2.3 and Cor. 6.4.1]{gh}, we know the first diagram is commutative. \\
\indent Now, since $d\log$ is surjective by definition, the map $\pi_{i}(\K_{\geq 0}(S))/p^{r}\to \pi_{i}(\K_{\geq 0}(S,\bbZ_{p}))/p^{r}$ is surjective and hence an isomorphism. In particular, each $\pi_{i-1}(\K_{\geq 0}(S))$ is $p$-torsion free as argued in the first paragraph, and the map $\pi_{i}(\K_{\geq 0}(S))/p^{r}\simeq\pi_{i}(\K_{\geq 0}(S)/p^{r})\to W_{r}\Omega^{i}_{S,\log}$ of question is an isomorphism. 
\end{proof}
\end{corollary}

\begin{remarkn}\label{rmk:padicmodprequiv}
Let $S$ be a local Cartier smooth ring. \\
(1) The $p$-torsion freeness part of Corollary \ref{cor:cartsmglmodpr} for each $\K_{i}(S)$ ($i\geq 0$) can also be shown directly; the following argument was provided to the author by Matthew Morrow. Consider the commutative diagram
\begin{equation*}
\begin{tikzcd}
& \K^{M}_{i}(S)/p \arrow{ld} \arrow{d} \arrow[bend left=60, two heads]{dd}{d\log} \\
\K_{i}(S)/p \arrow[hook]{r} & \K_{i}(S,\bbZ/p) \arrow{d}{\simeq} \\
& \Omega^{i}_{S,\log}.
\end{tikzcd}
\end{equation*}
Since $d\log$ is surjective by definition, the upper vertical arrow is surjective, and hence the horizontal map is surjective (equivalently an isomorphism). \\
(2) Hence, the mod $p^{r}$-versions imply the $p$-adic version. Since the pro-object $\{\pi_{i+1}(\K_{\geq 0}(S)/p^{r})\}_{r}$ is isomorphic to $\{\K_{i+1}(S)/p^{r}\}_{r}$ by (1), we in particular know its transition maps are surjective. Hence, the Milnor exact sequence gives an isomorphism $\pi_{i}\K_{\geq 0}(S,\bbZ_{p})\to\lim_{r}\pi_{i}(\K_{\geq 0}(S)/p^{r})$. Composing with the cyclotomic trace induced map, we obtain an isomorphism $\pi_{i}\K(S,\bbZ_{p})\to W\Omega^{i}_{S,\log}$. 
\end{remarkn}

\begin{remarkn}
Although \cite{km} proves Corollary \ref{cor:cartsmglmodpr} through a different method, it was indicated in \emph{loc. cit.} that an approach through the motivic filtration would be possible. 
\end{remarkn}

Let us record one application of Kelly-Morrow's generalization of the Geisser-Levine theorem for perfect rings.

\begin{lemma}\label{lem:perfdW}
Let $S$ be a perfect ring over $\bbF_{p}$. Then, there is a natural isomorphism $W\Omega^{\ast}_{S}\simeq W(S)$ of (strict) Dieudonn\'e algebras. 
\begin{proof}
Here, let us provide a proof using the formalism of saturated de Rham-Witt complexes. From $(L_{W(S)/\bbZ})/p\simeq S\dt_{W(S)}L_{W(S)/\bbZ}\simeq L_{S/\bbF_{p}}\simeq 0$, we know the $p$-completed cotangent complex $L_{W(S)/\bbZ}^{\wedge_{p}}$ vanishes. In particular, we know the $p$-completed de Rham complex $\widehat{\Omega}^{\ast}_{W(S)}$ is isomorphic to $W(S)$ as a Dieudonn\'e algebra. Note that $W(S)$ is a strict Dieudonn\'e algebra \cite[Ex. 2.5.6]{blm}. Now, for any strict Dieudonn\'e algebra $A^{\ast}$, we have natural isomorphisms $\Mor_{\CAlg_{\bbF_{p}}}(S,A^{0}/VA^{0})\simeq\Mor_{F}(W(S),A^{0})\simeq \Mor_{\text{DA}}(\widehat{\Omega}^{\ast}_{W(S)},A^{\ast})\simeq\Mor_{\text{DA}_{\text{str}}}(W(S),A^{\ast})$ by \cite[Prop. 3.6.3 and Var. 3.3.1]{blm}. Combined with Proposition \ref{prop:cartsmsaturated}, we know $\mathcal{W}\Omega^{\ast}_{S}\simeq W\Omega^{\ast}_{S}\simeq W(S)$ as strict Dieudonn\'e algebras.
\end{proof}
\end{lemma}

\begin{remarkn}
The following direct argument for a proof of Lemma \ref{lem:perfdW} was suggested to the author by Kay R\"ulling. First, observe that $\Omega^{i}_{W_{n}(S)}\simeq 0$ for each $i, n\geq1$.\footnote{For instance, one can use the fiber sequence $W_{n}(S)\dt_{\bbZ/p^{n}}L_{(\bbZ/p^{n})/\bbZ}\to L_{W_{n}(S)/\bbZ}\to L_{W_{n}(S)/(\bbZ/p^{n})}$ associated with the composition $\bbZ\to \bbZ/p^{n}\to W_{n}(S)$ and the vanishing $L_{W_{n}(S)/(\bbZ/p^{n})}\simeq L_{W(S)/\bbZ}^{\wedge_{p}}/p^{n}\simeq0$ as in the proof above.} On the other hand, by the construction of de Rham-Witt complexes, we have a natural surjection $\Omega^{i}_{W_{n}(S)}\twoheadrightarrow W_{n}\Omega^{i}_{S}$, cf. \cite[Th. I.1.3]{ill}, which forces the target to be zero. 
\end{remarkn}

\begin{proposition}\label{prop:Hiller}(Hiller)
Let $S$ be a perfect ring over $\bbF_{p}$. Then, $\K_{i}(S)$ is a $\bbZ[1/p]$-module for all $i>0$. 
\begin{proof}    
Since $S$ is a filtered colimit of colimit perfections of finitely presented perfect $\bbF_{p}$-algebras (as the colimit perfection is a left adjoint functor), we can assume that $S$ itself is of the form $S'_{\text{perf}}$ for some $S'$ finite type over $\bbF_{p}$. In particular, the underlying topological space of $\Spec S$ is the same as that of $\Spec S'$. By Zariski descent spectral sequence (note that $\K$ is a Postnikov complete sheaf by \cite[Th. 3.2]{etalek}), we are further reduced to the case of $S$ being local. By Corollary \ref{cor:cartsmglmodpr}, we know that each $\K_{i}(S)$ is $p$-torsion free and that $\K_{i}(S)/p\simeq \K_{i}(S,\bbZ_{p})/p\simeq W\Omega^{i}_{S,\log}/p$, which is zero by Lemma \ref{lem:perfdW}. Thus, $p$ acts invertibly on $\K_{i}(S)$.
\end{proof}
\end{proposition}

\begin{remark}
After the argument above was written down, the author was informed later that the idea of explaining Proposition \ref{prop:Hiller} through Geisser-Levine theorem appeared about two months earlier in \cite{amm}. It might be interesting to see if the approach using de Rham-Witt complexes can be generalized to cover certain mixed characteristic cases. 
\end{remark}

\section{More on prismatic cohomology of Cartier smooth rings}
\label{sec:cartsmprismetc}
In this section, we collect some consequences of \ref{sec:cartsmderW} relevant to prismatic cohomology complexes of Cartier smooth rings and their $p$-torsion free liftings following \cite{bl}. In the stated results, reference to \cite{blm} or \cite{bl} in their titles indicate that the corresponding results stated in there for regular Noetherian $\bbF_{p}$-algebras or ind-smooth maps remain valid for Cartier smooth rings or Cartier smooth maps; we explain how necessary modifications can be made.  \\
\indent First, note that Proposition \ref{prop:cartsmLdW} can be restated as follows:

\begin{proposition}[Proposition \ref{prop:cartsmLdW}]
Let $S$ be a Cartier smooth $\bbF_{p}$-algebra. Then, there is a natural equivalence $\Fil^{\sbullet}_{\N}\prism_{S}\simeq \Nyg^{\geq\sbullet}W\Omega_{S}$ in $\CAlg(\DFh(\bbZ_{p}))$. 
\begin{proof}
Recall that for $R\in\CAlg_{\bbF_{p}}^{\an}$, there is a natural equivalence $\Fil^{\sbullet}_{\N}\prism_{R}\simeq \Nyg^{\geq\sbullet}LW\Omega_{R}$; for instance, both sides commute with sifted colimits in $\CAlg^{\an}_{\bbF_{p}}$, and the equivalence can be checked on smooth $\bbF_{p}$-algebras. Specializing to the case $R=S$ and applying Proposition \ref{prop:cartsmLdW}, we have the equivalence $\Fil^{\sbullet}_{\N}\prism_{S}\simeq \Nyg^{\geq\sbullet}W\Omega_{S}$. 
\end{proof}
\end{proposition}

In particular, through the equivalence $\gamma^{\crys}_{\prism}:\prism_{S}\simeq\RG_{\crys}(\Spec S/\bbZ_{p})$ given by the crystalline comparison map \cite[Th. 4.6.1]{bl}, we have an equivalence $\RG_{\crys}(\Spec S/\bbZ_{p})\simeq W\Omega_{S}$, i.e., the derived crystalline cohomology of $S$ is computed by the de Rham-Witt complex of $S$. Also, recall that for Cartier smooth rings, the saturated de Rham-Witt complex can be identified with the (classical) de Rham-Witt complex via Proposition \ref{prop:cartsmsaturated}. Let us summarize what we have observed so far:

\begin{corollary}\label{cor:cartsmprism}
Let $S$ be a Cartier smooth $\bbF_{p}$-algebra. Then, there are natural equivalences
\begin{equation*}
\prism_{S}\simeq \RG_{\crys}(\Spec S/\bbZ_{p})\simeq LW\Omega_{S}\simeq \mathcal{W}\Omega_{S}\simeq W\Omega_{S}
\end{equation*}
in $\CAlg(\Dh(\bbZ_{p}))$. These equivalences respect filtrations, and lifts to equivalences
\begin{equation*}
\Fil^{\sbullet}_{\N}\prism_{S}\simeq \Nyg^{\geq\sbullet}LW\Omega_{S}\simeq \Nyg^{\geq\sbullet}\mathcal{W}\Omega_{S}\simeq\Nyg^{\geq\sbullet}W\Omega_{S}
\end{equation*} 
in $\CAlg(\DFh(\bbZ_{p}))$. 
\end{corollary}

By Zariski descent, one has analogous statements for Cartier smooth $\bbF_{p}$-schemes. We remark that the identification $\RG_{\crys}(\Spec S/\bbZ_{p})\simeq W\Omega_{S}$ can also be directly verified through the arguments provided in \cite{blm}:
 
\begin{proposition}\label{lem:cartsmcrysdW}(\cite[Th. 10.1.1]{blm} for smooth algebras over perfect fields) 
Let $S$ be a Cartier smooth $\bbF_{p}$-algebra. Then, there is a natural equivalence $\RG_{\crys}(\Spec S/\bbZ_{p})\simeq W\Omega_{S}$ in $\CAlg(\widehat{\mathcal{D}}(\bbZ_{p}))$. 
\begin{proof}
It suffices to check that there is a natural equivalence $\RG_{\crys}(\Spec S/\bbZ_{p})\simeq LW\Omega_{S}$ due to Proposition \ref{prop:cartsmLdW}, and this follows from \cite[Prop. 10.2.16 and Prop. 10.3.1]{blm}. For convenience, let us briefly explain their contents here. Write the restriction of $LW\Omega_{(-)}$ on quasisyntomic $\bbF_{p}$-algebras as $A(-)\in\CAlg(\Fun(\CAlg^{\QSyn}_{\bbF_{p}},\widehat{\mathcal{D}}(\bbZ_{p})))$. After the base change along $\widehat{\mathcal{D}}(\bbZ_{p})\to\mathcal{D}(\bbF_{p})$, it admits an equivalence $L\Omega_{(-)}\to\bbF_{p}\dt_{\bbZ_{p}}A(-)$ as commutative algebra objects of $\Fun(\CAlg^{\QSyn}_{\bbF_{p}},\mathcal{D}(\bbF_{p}))$. By \cite[Prop. 10.2.16]{blm}, there is a morphism $\RG_{\crys}(-/\bbZ_{p})\to A(-)$ in $\CAlg(\Fun(\CAlg^{\QSyn}_{\bbF_{p}},\widehat{\mathcal{D}}(\bbZ_{p})))$ which intertwines augmentation maps $\RG_{\crys}(-/\bbZ_{p})\xrightarrow{\epsilon^{\crys}}id$ and $A(-)\xrightarrow{\epsilon}id$. To show the morphism is an equivalence, it suffices to check its image in $\CAlg(\Fun(\CAlg^{\QSyn}_{\bbF_{p}},\mathcal{D}(\bbF_{p})))$ is an equivalence. This reduction mod $p$ gives a morphism $L\Omega_{(-)}\to L\Omega_{(-)}$ in $\CAlg(\Fun(\CAlg^{\QSyn}_{\bbF_{p}},\mathcal{D}(\bbF_{p})))$ compatible with the augmentation map $L\Omega_{(-)}\xrightarrow{\epsilon^{\dR}}id$. Now, \cite[Proof of Prop. 10.3.1]{blm} assures the morphism is an equivalence. 
\end{proof}
\end{proposition}

Let us remark on a few properties of $p$-torsion-free liftings of Cartier smooth $\bbF_{p}$-algebras. Note that these commutative rings are instances of F-smooth rings introduced and studied in \cite{bm}; see \emph{loc. cit.} for more general results hold for F-smooth rings. 

\begin{corollary}(cf. \cite[Rem. 5.4.3]{bl})
Let $R$ be a $p$-torsion-free commutative ring whose mod-$p$ quotient $\overline{R} := R/pR$ is a Cartier smooth $\bbF_{p}$-algebra. Then, there are natural equivalences $W\Omega_{\overline{R}}\simeq \prism_{\overline{R}}\simeq\dRh_{R}\simeq\widehat{\Omega}_{R}$ in $\CAlg(\Dh(\bbZ_{p}))$. In particular, for a $p$-adic formal scheme $\mathcal{X}$ flat over $\Spf\bbZ_{p}$ whose special fiber $X = \mathcal{X}\times_{\Spf\bbZ_{p}}\Spec\bbF_{p}$ is Cartier smooth, there is a natural equivalence $\RG_{\crys}(X/\bbZ_{p})\simeq\RG_{\dR}(\mathcal{X})$. 
\begin{proof}
Combine Corollary \ref{cor:cartsmprism}, absolute de Rham comparison theorem \cite[Th. 5.4.2]{bl}, and \cite[Prop. E.12]{bl}.
\end{proof}
\end{corollary}

\begin{remarkn}(cf. \cite[Rem. 4.7.4]{bl})\label{rmk:cartsmliftconj}
Let $R$ be a $p$-torsion-free commutative ring whose quotient $R/pR$ is Cartier smooth. Then, by \cite[Prop. E.12]{bl}, $L\widehat{\Omega}^{n}_{R}\simeq\widehat{\Omega}^{n}_{R}$ for all $n\geq0$. In particular, one can compute the conjugate filtration $\Fil^{\conj}_{\leq\sbullet}\widehat{\Omega}^{\slashed{D}}_{R}$ of the $p$-complete diffracted Hodge cohomology $R$; by the previous computation $\gr^{\conj}_{n}\widehat{\Omega}^{\slashed{D}}_{R}\simeq L\widehat{\Omega}^{n}_{R}[-n]$ is concentrated in degree $n$ for all $n\in\bbZ_{\geq0}$, and hence $\Fil^{\conj}_{\leq \sbullet}\widehat{\Omega}^{\slashed{D}}_{R}\simeq \tau^{\leq \sbullet}\widehat{\Omega}^{\slashed{D}}_{R}$.  
\end{remarkn}

The following is just a special case of a more general fact that F-smooth rings (by definition) have complete Nygaard filtrations, cf. \cite[Prop. 4.16]{bm}. 

\begin{proposition}(cf. \cite[Prop. 5.8.2]{bl})\label{prop:cartsmliftnygcpl}
Let $R$ be a $p$-torsion-free commutative ring whose quotient $R/pR$ is Cartier smooth. Then, for each $n\in\bbZ$, the prismatic complex $\prism_{R}\{n\}$ is Nygaard-complete. 
\begin{proof}
Let us explain how the proof \emph{loc. cit.} can be adjusted to prove the claimed case. The result follows from \cite[Lem. 5.8.3]{bl}, which in turn relies on \cite[Cor. 5.2.11]{bl}:
\begin{lemma}(cf. \cite[Cor. 5.2.11]{bl})
Let $(A,I)$ be a bounded prism and let $R$ be a $p$-completely flat $\overline{A}$(=$A/I$)-algebra for which $\overline{A}/p\overline{A}\to R/pR$ is a Cartier smooth map (see Lemma \ref{lem:cartsmflatbch} and above). Then, $F^{\ast}\prism_{R/A}$ is Nygaard-complete. 
\end{lemma}
Although the statement \emph{loc. cit.} is made for ind-smooth $\overline{A}/p\overline{A}$-algebras, the proof works without any change. As explained in the proof of \cite[Prop. 5.8.2 and Lem. 5.8.3]{bl}, the question is reduced to the statement that the relative Nygaard filtration on $F^{\ast}\prism_{(\overline{A}^{m}\otimes R)/A^{m}}$ is complete. By Lemma above, it suffices to check that for each $m\geq 0$, the natural map $\overline{A}^{m}/p\overline{A}^{m}\to (\overline{A}^{m}\otimes R)/p(\overline{A}^{m}\otimes R)$ is Cartier smooth.\footnote{Recall that $(A^{0},I^{0})$ is the transverse prism $(\bbZ[a_{0},a_{1}^{\pm1},a_{2},...]^{\wedge}_{(p,a_{0})}, (a_{0}))$ and $(A^{\sbullet},I^{\sbullet})$ is the associated cosimplicial (transverse) prism given by coproducts of copies of $(A^{0},I^{0})$, and $\overline{A}^{m} = A^{m}/I^{m}$.} From the $p$-torsion-freeness of $R$ and the fact that $(A^{m},I^{m})$ is transverse, we know $\overline{A}^{m}\otimes (R/pR)\simeq (\overline{A}^{m}\otimes R)/p(\overline{A}^{m}\otimes R)$, and hence $\overline{A}^{m}/p\overline{A}^{m}\otimes_{\bbF_{p}} (R/pR)\simeq (\overline{A}^{m}\otimes R)/p(\overline{A}^{m}\otimes R)$. In particular, by Lemma \ref{lem:cartsmflatbch}, the natural map $\overline{A}^{m}/p\overline{A}^{m}\to (\overline{A}^{m}\otimes R)/p(\overline{A}^{m}\otimes R)$ is Cartier smooth.  
\end{proof}
\end{proposition}

\begin{proposition}(cf. \cite[Prop. 5.7.9]{bl})
Let $R$ be a $p$-torsion-free commutative ring whose quotient $R/pR$ is Cartier smooth. Then, for each $n\in\bbZ$, the map $\Fil^{\sbullet}(\varphi\{n\}):\Fil^{\sbullet}_{\N}\prism_{R}\{n\}\to\prism_{R}^{[\sbullet+n]}$ exhibits the source as the connective cover of the target with respect to the Beilinson t-structure. 
\begin{proof}
The proof \emph{loc. cit.} works \emph{ad verbum}, using Remark \ref{rmk:cartsmliftconj} and Propsition \ref{prop:cartsmliftnygcpl} above. 
\end{proof}
\end{proposition}

Recall that a map $A\to B$ of commutative $\bbF_{p}$-algebras is Cartier smooth \cite[Def. E.10]{bl} if it is flat, $L_{B/A}$ is a flat discrete $B$-module, and the inverse Cartier map $\Omega^{i}_{B^{(1)}/A}\to H^{i}(\Omega^{\ast}_{B/A})$ is an isomorphism of $B^{(1)} = A\otimes_{\varphi,A}B$-modules for all $i\geq 0$. 

\begin{lemma}\label{lem:cartsmflatbch}
Let $A\to B$ a map of commutative $\bbF_{p}$-algebras and let $A'\to B'$ be its base change along a flat map $A\to A'$. If $A\to B$ is Cartier smooth, then so is $A'\to B'$. If $A\to A'$ is faithfully flat, then the converse holds. 
\begin{proof}
Since properties of being flat and of being an isomorphism are preserved under base change and satisfy faithfully flat descent, it suffices to check that $B'^{(1)}\otimes_{B^{(1)}}C^{-1}_{B/A}:B'^{(1)}\otimes_{B^{(1)}}\Omega^{i}_{B^{(1)}/A}\to B'^{(1)}\otimes_{B^{(1)}}H^{i}(\Omega^{\ast}_{B/A})$ agrees with  $C^{-1}_{B'/A'}:\Omega^{i}_{B'^{(1)}/A'}\to H^{i}(\Omega^{\ast}_{B'/A'})$ as graded $B'^{(1)}$-algebra maps (where $C^{-1}$ denotes the inverse Cartier map). Consider the commutative diagram
\begin{equation*}
\begin{tikzcd}
A\arrow{rrr}{\varphi}\arrow{dd}\arrow{rd}&   &           &    A \arrow{rd} \arrow{dd} \arrow{ld}& \\
 & B \arrow{r} &  B^{(1)} \arrow[rr, crossing over]&      & B \arrow{dd}\\
A'\arrow{rrr}{\varphi} \arrow{rd}&  &            &   A' \arrow{rd}\arrow{ld} & \\
 & B'\arrow{r} \arrow[from=uu, crossing over]&  B'^{(1)} \arrow{rr} \arrow[from=uu, crossing over]&      & B',
\end{tikzcd}
\end{equation*}
where the four vertical edges of the outer cube are Frobenius. By definition the left, upper, and bottom sides of the left hexahedron are pushout squares, and hence the right side of it is a pushout square. In particular, we have $B'^{(1)}\otimes_{B^{(1)}}\Omega^{i}_{B^{(1)}/A}\simeq\Omega^{i}_{B'^{(1)}/A'}$. As the right side of the outer cube is a pushout square, we know the front square of the right prism is a pushout square. Combined with $\Omega^{\ast}_{B/A}\simeq\Omega^{\ast}_{B/B^{(1)}}$ and $\Omega^{\ast}_{B'/A'}\simeq\Omega^{\ast}_{B'/B'^{(1)}}$ \cite[tag 0CCC]{stacks}, this gives $B'^{(1)}\otimes_{B^{(1)}}\Omega^{\ast}_{B/A}\simeq \Omega^{\ast}_{B'/A'}$. These identify source and target of both of the maps (using the flatness assumption). To prove the claim, it suffices to check that under these identifications, $B'^{(1)}\otimes_{B^{(1)}}C^{-1}_{B/A}$ agrees with $C^{-1}_{B'/A'}$ on degree 1 elements of the form $d(1\otimes g')$ for each $g'\in B'$. Let $g' = \sum_{j}f'_{j}\otimes g_{j}\in B' = A'\otimes_{A}B$. Its image in $B'^{(1)}$ is $1\otimes g' = \sum_{j}1\otimes (f'_{j}\otimes g_{j}) = \sum_{j}\varphi(f'_{j})\cdot 1\otimes (1\otimes g_{j})$, and hence $d(1\otimes g') = \sum_{j}\varphi(f'_{j})d(1\otimes (1\otimes g_{j}))$. The image of the latter element by $C^{-1}_{B'/A'}$ is $\sum_{j}\varphi(f'_{j})\cdot[(1\otimes g_{j})^{p-1}d(1\otimes g_{j})]$. On the other hand, $d(1\otimes g')$ corresponds to $\sum_{j}\varphi(f'_{j})\otimes d(1\otimes g_{j})$ in $B'^{(1)}\otimes_{B^{(1)}}\Omega^{1}_{B^{(1)}/A}$, and its image by $B'^{(1)}\otimes_{B^{(1)}}C^{-1}_{B/A}$ is $\sum_{j}\varphi(f'_{j})\otimes [g_{j}^{p-1}dg_{j}]$, which in turn corresponds to $C^{-1}_{B'/A'}(d(1\otimes g'))$ as computed above.  
\end{proof}
\end{lemma}

\appendix

\section{Filtrations}
\label{sec:filt}

We review the notion of filtered objects and study the standard spectral sequence associated with filtered objects in Proposition \ref{prop:filtspsq}. 

\begin{definition}
Let $\mathcal{D}$ be a stable presentable $\infty$-category.\\
(1) A \emph{filtered object} in $\mathcal{D}$ is a $\mathbb{Z}^{\op}$-shaped diagram $F^{\geq\sbullet}:\N\bbZ^{\op}\to\mathcal{D}$ in $\mathcal{D}$, where $\bbZ$ as a poset is equipped with the usual partial order. We denote $\mathcal{D}^{\fil}:=\Fun(\bbZ^{\op},\mathcal{D})$ for the $\infty$-category of filtered objects in $\mathcal{D}$. \\
(2) A \emph{graded object} in $\mathcal{D}$ is a $\bbZ^{\disc}$-shaped diagram $X:\N\bbZ^{\disc}\to\mathcal{D}$ in $\mathcal{D}$, where $\bbZ^{\disc}$ is the discrete category whose underlying set of objects is $\bbZ$. We denote $\mathcal{D}^{\gr}:=\Fun(\bbZ^{\disc},\mathcal{D})$ for the $\infty$-category of graded objects in $\mathcal{D}$. 
\end{definition}

In particular, we will only consider $\bbZ$-graded objects and (descending or increasing) $\bbZ$-filtered objects in this article. For notational convenience, let us deal with the case of descending filtrations in this section. Consider the unique map $p:\N\bbZ^{\op}\to\Delta^{0}$ of simplicial sets and the induced functor $p^{\ast}:\mathcal{D}\to\mathcal{D}^{\fil}$ which assigns each object of $\mathcal{D}$ to a constant filtration. As $p$ is a Cartesian and coCartesian fibration, we have adjunctions $p\lsh\dashv p^{\ast}\dashv p_{\ast}$ given by the functors of left and right Kan extensions along $p$, which are computed as $p\lsh(F)\simeq\colim_{n\to-\infty}F^{\geq n}$ and $p_{\ast}(F)\simeq \lim_{n\to\infty}F^{\geq n}$ respectively.

\begin{definition}
Let $F^{\geq\sbullet}$ be a filtered object in a stable presentable $\infty$-category $\mathcal{D}$. \\
(1) $F^{\geq\sbullet}$ is \emph{complete} if $p_{\ast}(F)\simeq 0$, i.e., $\lim_{n\to\infty}F^{\geq n}\simeq0$. \\
(2) An object $F\in\mathcal{D}$ is an \emph{underlying object} of $F^{\geq\sbullet}$ if there is an equivalence $F\simeq p\lsh(F)$, i.e., $\colim_{n\to-\infty}F^{\geq n}\simeq F$. In this case, $F^{\geq\sbullet}$ is called an \emph{exhausive filtration} of $F$.  
\end{definition} 

Above definitions formalize the idea of linear algebraic filtrations prevalent in algebra and geometry. An important operation useful in studying filtrations is taking the associated graded object of a given filtered object, which we will review below shortly; this often makes the study of the underlying object more accessible by providing further linear algebraic data or a more commutative object. Conceptually, this operation also provides a deformation picture, that each $\bbZ$-filtered object gives a one-parameter deformation of its associated graded object (i.e., the value at the special fiber) by its underlying object (i.e., the value at the generic fiber). 

\begin{example}
An archetypal example is the sheaf of differential operators $\mathcal{D}_{X}$ on a smooth variety $X$ over a characteristic 0 field $k$, equipped with the natural exhaustive increasing $\bbZ$-filtration $\Fil_{\leq\sbullet}\mathcal{D}_{X}$ by order. The associated graded object is naturally isomorphic to $\Sym^{\sbullet} T_{X/k}$ as a (graded Poisson) $\mathscr{O}_{X}$-algebra, and hence $\mathcal{D}_{X}$ is a one-parameter deformation or a \emph{quantization} of the cotangent bundle of $X$; this is important in algebraic geometry and representation theory. Analogous statements can be made for their modules, e.g., for modules whose underlying sheaves are vector bundles on $X$, flat connections equipped with a filtration compatible with the filtration on $\mathcal{D}_{X}$ deform graded Higgs bundles. 
\end{example}

\begin{definition}
Let $F^{\geq\sbullet}\in\mathcal{D}^{\fil}$. Its \emph{associated graded} object $\gr(F)\in\mathcal{D}^{\gr}$ is defined as $\gr(F)^{n} = \gr^{n}(F):=\cof(F^{\geq n+1}\to F^{\geq n})$ for each $n\in\bbZ$. This construction naturally defines a functor $\gr:\mathcal{D}^{\fil}\to\mathcal{D}^{\gr}$.  
\end{definition}

By construction $\gr$ commutes with small colimits, and hence is a left adjoint functor under our assumption on $\mathcal{D}$. Its right adjoint can be described as follows:

\begin{proposition}\label{prop:grtadj}
Let $\mathcal{D}$ be a stable presentable $\infty$-category. Consider the functor $t:\mathcal{D}^{\gr}\to\mathcal{D}^{\fil}$ naturally described as $X\mapsto (\cdots\xrightarrow{0} X^{n+1}\xrightarrow{0} X^{n}\xrightarrow{0} X^{n-1}\xrightarrow{0}\cdots)$, i.e., $t(X)^{\geq n}\simeq X^{n}$ and $(t(X)^{\geq n+1}\to t(X)^{\geq n})\simeq 0$ for all $n\in\bbZ$. Then, $t$ is a right adjoint functor, and it fits into an adjoint situation $\gr\dashv t:\mathcal{D}^{\gr}\to\mathcal{D}^{\fil}$. 
\begin{proof}
First, for each $n\in\bbZ$, consider the functor $c_{n}:\mathcal{D}\to\mathcal{D}^{\fil}$ naturally described by sending $Y\in\mathcal{D}$ to $(\cdots\to 0\to Y\to 0\to\cdots)$, i.e., $c_{n}(Y)^{\geq n}\simeq Y$ and $c_{n}(Y)^{\geq i}\simeq 0$ for $i\neq n$.
\begin{lemma}\label{lem:grnadj}
There is an adjoint situation $\gr^{n}\dashv c_{n}:\mathcal{D}\to\mathcal{D}^{\fil}$.
\begin{proof}
It suffices to check that there is a natural equivalence $\Map_{\mathcal{D}^{\fil}}(F,c_{n}(Y))\simeq\Map_{\mathcal{D}}(\gr^{n}(F),Y)$ for $F\in\mathcal{D}^{\fil}$ and $Y\in\mathcal{D}$. Recall that mapping spaces of $\mathcal{D}^{\fil} = \Fun(\N\mathbb{Z}^{\op},\mathcal{D})$ can be described through an end construction \cite[Prop. 2.3]{gla}. In our case, the twisted arrow category $\mathrm{Tw}(\bbZ^{\op})$ is given by (the nerve of) $\{(p,q)~|~p\geq q\}\subseteq\bbZ\times\bbZ^{\op}$. Hence, 
\begin{align*} 
\Map_{\mathcal{D}^{\fil}}(F,c_{n}(Y)) &\simeq \int_{\bbZ^{\op}}\Map_{\mathcal{D}}(F^{\geq\sbullet},c_{n}(Y)^{\geq\sbullet})\simeq\lim_{\{(p,q)\in\bbZ\times\bbZ^{\op}~|~p\geq q\}}\Map_{\mathcal{D}}(F^{\geq p}, c_{n}(Y)^{\geq q})\\
&\simeq \fib\left(\Map_{\mathcal{D}}(F^{\geq n},c_{n}(Y)^{\geq n})\to \Map_{\mathcal{D}}(F^{\geq n+1},c_{n}(Y)^{\geq n})\right)\\
&\simeq \Map_{\mathcal{D}}(\cof(F^{\geq n+1}\to F^{\geq n}),c_{n}(Y)^{\geq n})\simeq \Map_{\mathcal{D}}(\gr^{n}(F), Y). 
\end{align*}
While transitioning to the second line, we used that $c_{n}(Y)^{\geq q}\simeq 0$ for $q\neq n$.
\end{proof}
\end{lemma}
Now, we can check that there is a natural equivalence $\Map_{\mathcal{D}^{\gr}}(\gr F, X)\simeq \Map_{\mathcal{D}^{\fil}}(F, t(X))$ for $F = F^{\geq\sbullet}\in\mathcal{D}^{\fil}$ and $X\in\mathcal{D}^{\gr}$. In fact, from $\mathcal{D}^{\gr}\simeq\prod_{\bbZ}\mathcal{D}$, we have 
\begin{align*}
\Map_{\mathcal{D}^{\gr}}(\gr F, X)&\simeq \prod_{n\in \bbZ}\Map_{\mathcal{D}}(\gr^{n}(F), X^{n})\\
&\simeq \prod_{n\in\bbZ} \Map_{\mathcal{D}^{\fil}}(F, c_{n}(X^{n}))\simeq \Map_{\mathcal{D}^{\fil}}(F,\prod_{n\in\bbZ}c_{n}(X^{n}))\simeq \Map_{\mathcal{D}^{\fil}}(F,t(X)).
\end{align*}
Here, for the second equivalence we used the adjunction of Lemma \ref{lem:grnadj}, and for the final equivalence we used the natural equivalence $t(X)\simeq \prod_{n\in\bbZ}c_{n}(X^{n})$ in $\mathcal{D}^{\fil}$. 
\end{proof}
\end{proposition}

\begin{remarkn}
Suppose $\mathcal{D}$ is a symmetric monoidal stable presentable $\infty$-category. Then, $\mathcal{D}^{\fil}$ and $\mathcal{D}^{\gr}$ admit symmetric monoidal structures via Day convolution operations\footnote{Here, both $\N\bbZ^{\op}$ and $\N\bbZ^{\disc}$ are regarded as symmetric monoidal $\infty$-categories via addtion $+$ on $\bbZ$.}. Moreover, arguments of \cite[Sec. 3.2]{wald}, applied \emph{mutatis mutandis}, show that $\gr$ is symmetric monoidal, and hence the adjunction of Proposition \ref{prop:grtadj} enhances to a symmetric monoidal adjunction. 
\end{remarkn}

\begin{proposition}\label{prop:filtspsq}
Let $\mathcal{D}$ be a stable presentable $\infty$-category equipped with a t-structure compatible with filtered colimits (i.e., $\mathcal{D}_{\geq 0}$ is Grothendieck prestable). Suppose $F^{\geq\sbullet}$ is a filtered object of $\mathcal{D}$ which is complete, and suppose $F\in\mathcal{D}$ is its underlying object (i.e., $F^{\geq\sbullet}$ is an exhaustive filtration on $F$). Then, there is a conditionally convergent spectral sequence $E^{i,j}_{2} = \pi_{-i-j}(\gr^{-j}F^{\geq \sbullet})\Rightarrow \pi_{-i-j}F$ in $\mathcal{D}^{\heartsuit}$ (of cohomological grading). 
\begin{proof}
This is a very standard result used in a vast number of literatures in stable homotopy theory, often with slightly different indexings depending on the context. Nevertheless, let us state and prove the proposition here to clarify our convention. The spectral sequence can be constructed by considering the object $F(\sbullet):=F^{\geq -\sbullet}\in\Fun(\bbZ,\mathcal{D})$, and applying the construction of \cite[1.2.2.6]{ha} on $F(\sbullet)$ (by viewing it as a gapped object on $\bbZ\cup-\infty$ and then restricting to $\bbZ$). Equivalently, one can appeal to \cite{boardman}. In fact, fiber sequences $F^{\geq i+1}\to F^{\geq i}\to \gr^{i}F^{\geq\sbullet}$ induce an exact couple 
\begin{equation*}
\begin{tikzcd}[row sep=3mm, column sep=5mm]
\cdots\to\pi_{-\ast}F^{\geq i+1} \arrow{rr} & & \pi_{-\ast}F^{\geq i}\to\cdots \arrow{ld} \\
  & \pi_{-\ast}\gr^{i}F^{\geq\ast} \arrow{lu}{+1} &
\end{tikzcd}
\end{equation*}
in $\mathcal{D}^{\heartsuit}$, and the first page of the associated spectral sequence takes the form $E^{i,j}_{1} = \pi_{-i-j}\gr^{i}F^{\geq \sbullet}\xrightarrow{d_{1}} E^{i+1,j}_{1} = \pi_{-i-j-1}\gr^{i+1}F^{\geq\sbullet}$. To get the desired form, set $\ell=-i$ and $k = i+j-\ell = 2i+j$ (hence $i=-\ell$ and $j=k+2\ell$). As $(i,j)\mapsto (i+r,j-r+1)$ sends $(k,\ell)$ to $(k+r+1,\ell-r)$, we can define the new spectral sequence $((E')^{k,\ell}_{r},d_{r})_{r\geq 2}$ by $(E^{i,j}_{r-1} = E^{-\ell, k+2\ell}_{r-1}, d_{r-1})_{r\geq 2}$. In particular, the initial second page consists of objects $(E')^{k,\ell}_{2} = \pi_{-k-\ell}\gr^{-\ell}F^{\geq\sbullet}$. Since $\lim_{i\to\infty}F^{\geq i}\simeq 0$, we know $\lim^{p}_{i\to\infty}\pi_{-\ast}F^{\geq i}=0$ for $p=0,1$. Thus, our spectral sequence converges conditionally to the colimit $\pi_{-\ast}F$. 
\end{proof}
\end{proposition}

\printbibliography

\end{document}